\documentclass[reqno,11pt]{amsart}  
\usepackage{amsmath} 
\usepackage{amssymb} 
\usepackage{amsthm} 
\newcommand\NoBlackBoxes{\global\overfullrule0pt} 
\NoBlackBoxes 
\theoremstyle{plain} 
 
\def\4{\kern1pt}

\def\6{\vphantom0}

\def\8{\kern-10pt} 
\def\7#1{_{(#1)}}

\makeatletter 
\let\serieslogo@\relax 
\let\@setcopyright\relax 
 
\def\speciallabelmark#1{\def\@currentlabel{#1}} 
\makeatother 
 
\begin{document} 
 
\def\ffrac#1#2{\raise.5pt\hbox{\small$\4\displaystyle\frac{\,#1\,}{\,#2\,}\4$}} 
\def\ovln#1{\,{\overline{\!#1}}} 
\def\ve{\varepsilon} 
\def\kar{\beta_r} 

\title{NON-UNIFORM BOUNDS IN LOCAL LIMIT THEOREMS \\  
IN CASE OF FRACTIONAL MOMENTS 
} 
 
\author{S. G. Bobkov$^{1,4}$} 
\thanks{1) School of Mathematics, University of Minnesota, USA; 
Email: bobkov@math.umn.edu} 
\address 
{Sergey G. Bobkov \newline 
School of Mathematics, University of Minnesota  \newline  
127 Vincent Hall, 206 Church St. S.E., Minneapolis, MN 55455 USA 
\smallskip} 
\email {bobkov@math.umn.edu}  
 
\author{G. P. Chistyakov$^{2,4}$} 
\thanks{2) Faculty of Mathematics, University of Bielefeld, Germany; 
Email: chistyak@math.uni-bielefeld.de} 
\address 
{Gennady P. Chistyakov\newline 
Fakult\"at f\"ur Mathematik, Universit\"at Bielefeld\newline 
Postfach 100131, 33501 Bielefeld, Germany} 
\email {chistyak@math.uni-bielefeld.de} 
 
\author{F. G\"otze$^{3,4}$} 
\thanks{3) Faculty of Mathematics, University of Bielefeld, Germany; 
Email: goetze@math.uni-bielefeld.de} 
\thanks{4) Research partially supported by  
SFB 701} 
\address 
{Friedrich G\"otze\newline 
Fakult\"at f\"ur Mathematik, Universit\"at Bielefeld\newline 
Postfach 100131, 33501 Bielefeld, Germany} 
\email {goetze@mathematik.uni-bielefeld.de}


 
\begin{abstract} 
Edgeworth-type expansions for convolutions of probability densities 
and powers of the characteristic functions with non-uniform error terms 
are established for i.i.d. random variables with finite (fractional) moments  
of order $s \geq 2$, where $s$ may be noninteger. 
\end{abstract} 

\maketitle 
\markboth{S. G. Bobkov, G. P. Chistyakov and F. G\"otze}{Edgeworth-type expansions  
for characteristic functions}

\def\theequation{\thesection.\arabic{equation}} 
\def\E{{\bf E}} 
\def\R{{\bf R}} 
\def\C{{\bf C}} 
\def\P{{\bf P}} 
\def\H{{\rm H}} 
\def\Im{{\rm Im}} 
\def\Tr{{\rm Tr}} 
 
\def\k{{\kappa}} 
\def\M{{\cal M}} 
\def\Var{{\rm Var}} 
\def\Ent{{\rm Ent}} 
\def\O{{\rm Osc}_\mu} 
 
\def\ep{\varepsilon} 
\def\phi{\varphi} 
\def\F{{\cal F}} 
\def\L{{\cal L}} 
 
\def\be{\begin{equation}} 
\def\en{\end{equation}} 
\def\bee{\begin{eqnarray*}} 
\def\ene{\end{eqnarray*}}

 
\section{{\bf Introduction}} 
\setcounter{equation}{0} 
 
\vskip2mm 
Let $(X_n)_{n \geq 1}$ be independent and identically distributed random 
viariables with $\E X_1 = 0$ and $\E X_1^2 = 1$.  
 
If $X_1$ has finite moments of all orders, and if the densities $\rho_n$ of  
the normalized sums $S_n = (X_1 + \dots + X_n)/\sqrt{n}$ exist, they admit  
a formal Edgeworth-type expansion with respect to the powers of $1/\sqrt{n}$ 
\be 
\rho_n(x) = \varphi(x) + \sum_{k=1}^\infty q_k(x)\,n^{-k/2}. 
\en 
Here, $\varphi(x) = \frac{1}{\sqrt{2\pi}}\,e^{-x^2/2}$ denotes the density  
of the standard normal law, 
\be 
q_k(x) \ = \, \varphi(x)\, \sum H_{k + 2j}(x) \, 
\frac{1}{p_1!\dots p_k!}\, \bigg(\frac{\gamma_3}{3!}\bigg)^{p_1} \dots 
\bigg(\frac{\gamma_{k+2}}{(k+2)!}\bigg)^{p_k}, 
\en 
where $H_k(x)$ are the Chebyshev-Hermite polynomials, and  
$$ 
\gamma_k = i^{-k} \frac{d^k}{dt^k}\, \log \E\, e^{itX_1}\big|_{t=0} 
$$ 
denote the cumulants of the underlying distribution. The summation in (1.2)  
runs over all non-negative integer solutions $(p_1,\dots,p_k)$ to the equation 
$p_1 + 2 p_2 + \dots + k p_k = k$ with $j = p_1 + \dots + p_k$. 
 
A precise asymptotic statement about the formal series (1.1) requires that 
some moment $\E |X_1|^s$ of order $s \geq 2$ is finite 
(while the moments of higher orders may be infinite). In this case, the $k$-th order  
cumulants are well-defined for the values $k = 1,\dots,m$, and respectively,  
the functions $q_k$ are defined for $k \leq m-2$, where $m = [s]$ is the integer  
part of $s$. Therefore, one needs to evaluate the error of the approximation of  
$\rho_n$ by the following partial sums of the series (1.1), 
$$ 
\varphi_m(x) = \varphi(x) + \sum_{k=1}^{m-2} q_k(x)\,n^{-k/2}, \qquad m = [s]. 
$$ 
 
One of the aims of this paper is to prove the following: 
 
\vskip5mm 
{\bf Theorem 1.1.} {\it Assume that $\E\, |X_1|^s < +\infty$, for some $s \geq 2$. 
Suppose $S_{n_0}$ has a bounded density $\rho_{n_0}$ for some $n_0$. Then for all $n$ large 
enough, $S_n$ have continuous densities $\rho_n$ satisfying, as $n \rightarrow \infty$, 
\be 
(1 + |x|^m)\, \big(\rho_n(x) - \varphi_m(x)\big) = o\big(n^{-(s-2)/2}\big) 
\en 
uniformly for all $x$. Moreover, 
\be 
(1 + |x|^s) \big(\rho_n(x) - \varphi_m(x)\big) = o\big(n^{-(s-2)/2}\big) + 
\left(1+|x|^{s-m}\right)\,\big(O(n^{-(m-1)/2})+ o(n^{-(s-2)})\big). 
\en 
} 
 
\vskip2mm 
In fact, the implied sequence and constants in the error terms hold uniformly  
in the class of all densities with precribed moment tail function  
$t \rightarrow \E\, |X_1|^s\,1_{\{|X_1| > t\}}$, parameter $n_0$ and a bound on 
the density $\rho_{n_0}$. 
 
For $|x|$ of order 1, or when $s = m$ is integer, both relations are equivalent.  
But for large values of $|x|$ and $s>m$, the assertion (1.4) gives an  
improvement over (1.3), which is essential in some applications. 
 
If $2 \leq s < 3$, (1.4) becomes 
$$ 
(1 + |x|^s) \big(\rho_n(x) - \varphi(x)\big) = o\big(n^{-(s-2)/2}\big) + 
(1+|x|^{s-2})\, o\big(n^{-(s-2)}\big). 
$$ 
In particular, for the smallest value $s=2$, this contains the Gnedenko 
local limit theorem 
$\sup_x |\rho_n(x) - \varphi(x)| \rightarrow 0$, as $n \rightarrow \infty$. 
 
\vskip2mm 
If $s = m$ is integer and $m \geq 3$, Theorem 1.1 is well-known; (1.3)-(1.4) 
then simplify to 
\be 
(1 + |x|^m) \big(\rho_n(x) - \varphi_m(x)\big) = o\big(n^{-(m-2)/2}\big). 
\en 
In this formulation the result is due to Petrov [Pe1] (cf. also Petrov [Pe2],  
p.\,211, or Bhattacharya and Ranga Rao [B-RR], p.\,192). Without the term $1 + |x|^m$, 
(1.5) can be found in the classical book [G-K]; this weaker variant goes back  
to the results by Gnedenko [G] and an earlier work by Cramer  
(who used, according to [G-K], additional assumptions on the underlying density). 
 
Thus, Theorem 1.1 extends these well-known results to the case, where $s$ is not  
necessarily integer. The range $2 \leq s < 3$ is of interest, as well.  
Our interest in these somewhat technical extensions, especially (1.4),  
was motivated by open questions as to the actual rate of convergence  
in the so-called entropic central limit theorem. Here the relation (1.4) led  
to an unexpected behaviour of the error in the approximation of 
the entropy of sums of independent summands when $s$ increases from 2 to 4.  
(This error stabilizes at $s=4$, in contrast with the usual  
Berry-Esseen-type theorem for distribution functions, where stabilization of errors 
starts at $s=3$). In the entropic central limit theorem the classical non-uniform 
bound (1.5) is not precise enough to derive upper bounds for errors. 
 
Note that the assumption of boundedness of $\rho_n$ in Theorem 1.1  
(for some $n$ or, equivalently, for all large $n$) is necessary for  
conclusions such as (1.3)-(1.5). It is equivalent to the property  
that the characteristic function $v(t) = \E\, e^{itX_1}$ is integrable  
with some power $\nu \geq 1$, i.e., 
\be 
\int_{-\infty}^{+\infty} |v(t)|^\nu\,dt < +\infty. 
\en 
In this case $\rho_n$ are bounded for all $n \geq 2\nu$. The condition (1.6) 
is sometimes called "smoothness"; it appears naturally in many problems 
of the asymptotic behaviour of the densities (see e.g. [S] for detailed 
discussion). 
 
Nevertheless, this condition may be removed at all, if we require that  
(1.3)-(1.4) hold true for slightly modified densities, rather than for $\rho_n$. 
 
\vskip5mm 
{\bf Theorem 1.2.} {\it Let $\E\, |X_1|^s < +\infty$, for some $s \geq 2$. Let $c$ 
denote an arbitrary number with $0 < c < 1$. Suppose that for $n$ large enough $S_n$ have  
absolutely continuous distributions with densities  
$\rho_n$. Then, for some probability densities $\widetilde \rho_n$, 
 
\vskip2mm 
$a)$\, The relations $(1.3)$-$(1.4)$ hold true for $\widetilde \rho_n$  
in place of $\rho_n$ \!$;$ 
 
\vskip1mm 
$b)$\, $\int_{-\infty}^{+\infty} |\widetilde \rho_n(x) - \rho_n(x)|\,dx < c^n$, for all 
$n$ large enough \!$;$ 
 
\vskip1mm 
$c)$\, $\widetilde \rho_n(x) = \rho_n(x)$ almost everywhere, if $\rho_n$ is bounded 
$({\it a.e.})$ 
}

\vskip5mm 
It seems that Theorem 1.2 has not been stated in the literature,  
even when $s$ is integer. Here, the property $c)$ is added to include Theorem 1.1  
in Theorem~1.2 as a particular case.  
 
It turns out that the statement of Theorem 1.2 is more appropriate for a number of  
applications. For example, it implies that $\rho_n - \varphi_m \rightarrow 0$ in the mean,  
i.e., there is convergence in total variation norm for the corresponding distributions  
with rate 
\be 
\int_{-\infty}^{+\infty} |\rho_n(x) - \varphi_m(x)|\,dx  = o\big(n^{-(s-2)/2}\big). 
\en 
For $s=2$ and $\varphi_2(x) = \varphi(x)$, this statement corresponds  
to a theorem of Prokhorov [Pr], while for $s=3$ and  
$\varphi_3(x) = \varphi(x) \big(1 + \alpha_3\, \frac{x^3 - 3x}{6\sqrt{n}}\big)$  
-- to the result of Sirazhdinov and Mamatov [S-M] 
(they also covered the case $2 < s < 3$ with $O$ in place of $o$ in (1.7)). 
If $s \geq 3$ is integer, (1.7) is mentioned in [Pe2] for a more general $L^p$-convergence  
-- however, under the assumption that the densities $\rho_n$ are bounded. 
 
Theorem 1.2 allows to study non-uniform convergence in (1.3)-(1.4) as well when 
excluding exceptional "small" sets (via additional assumptions of entropic type). 
 
\vskip2mm 
The proof of Theorems 1.1 and 1.2, which are formally based on the application  
of the inverse Fourier transforms, involves operators, namely, Liouville fractional  
integrals and derivatives. 
For this step, we analyse the decay of the Fourier transform for special  
classes of finite measures with finite fractional moments. 
(Apparently standard truncation methods are much to density-sensitive and do  
not provide the required asymptotics.) 
An essential part of the argument is devoted to the routine analysis of powers  
of the characteristic functions and more general Fourier transforms 
in Edgeworth-type expansions. For this step the requirement (1.6) is irrelevant. 
 
In order to 
describe one of the main intermediate results, which is, as we believe, of an  
independent interest, let us start with a random variable $X$, such that  
$\E X = 0$, $\E X^2 = 1$, and $\E\, |X|^s < +\infty$, for some $s \geq 2$.  
Introduce the characteristic function $v(t) = \E\, e^{itX}$, $t \in \R$. 
 
If $m = [s]$, the normalized powers $v_n(t) = v(\frac{t}{\sqrt{n}})^n$, 
that is, the characteristic functions of $S_n$, can be approximated by the functions 
$$ 
u_m(t) = e^{-t^2/2}\,  
\bigg(1 + \sum_{k=1}^{m-2} P_k(it)\,n^{-k/2}\bigg). 
$$ 
Here we use the classical polynomials 
$$ 
P_k(t) \ \ = \sum_{p_1 + 2 p_2 + \dots + k p_k = k}   
\frac{1}{p_1!\dots p_k!}\, 
\bigg(\frac{\gamma_3}{3!}\bigg)^{p_1} \dots 
\bigg(\frac{\gamma_{k+2}}{(k+2)!}\bigg)^{p_k} 
t^{k + 2(p_1 + \dots + p_k)} 
$$ 
of degree $3k$, where the summation is performed as in (1.2). Another way to introduce 
these polynomials is to require that every $q_k(x)$ has the Fourier transform 
$e^{-t^2/2} P_k(it)$, so that $u_m(t)$ appears as the Fourier transform of $\varphi_m(x)$. 
 
The following statement is standard: in the interval $|t| \leq n^{1/6}$ 
$$ 
|v_n(t) - u_m(t)|\, \leq\, \frac{\ep_n}{n^{(m-2)/2}}\, 
\big(|t|^{m'} + |t|^{m''}\big)\,  e^{-t^2/2}, 
$$ 
where $\ep_n$ do not depend on $t$ and satisfy $\ep_n \rightarrow 0$, as  
$n \rightarrow \infty$ (with certain powers $m'$ and $m''$). 
Similar bounds also hold for the derivatives of orders $p = 1,\dots,m$, 
namely 
$$ 
\bigg|\frac{d^p}{dt^p}\, v_n(t) - \frac{d^p}{dt^p}\, u_m(t)\bigg|\, \leq\,  
\frac{\ep_n}{n^{(m-2)/2}}\,\big(|t|^{m'} + |t|^{m''}\big)\,  e^{-t^2/2}. 
$$ 
This bound is proved in Petrov [Pe1], cf. also [Pe2], pp. 209--211 
(for $m \geq 3$ and $|t| \leq n^{1/7}$). We refine this result with general  
values of $s \geq 2$ by proving that 
\be 
\bigg|\frac{d^p}{dt^p}\, v_n(t) - \frac{d^p}{dt^p}\ u_m(t) \bigg| \leq  
\frac{\ep_n}{n^{(s-2)/2}}\, \big(|t|^{m'} + |t|^{m''}\big)\, e^{-t^2/2}. 
\en 
However, the error term in this approximation is still not sufficiently  
small for our applications, and we have to look for other related  
representations of $v_n$. In analogue with $u_m$, introduce 
$$ 
e_m(t) = e^{-t^2/2}\, \bigg(1 + \sum_{k=1}^{m-2} P_k(it)\bigg). 
$$ 
 
\vskip2mm 
{\bf Theorem 1.3.} {\it Let $\E\, |X|^s < +\infty$ $(s \geq 2)$.  
For all $p = 0,1\dots,m$, and all $|t| \leq cn^{1/6}$, 
\be 
\frac{d^p}{dt^p}\, (v_n(t) - u_m(t)) = n\, \frac{d^p}{dt^p} 
\left[\bigg(v\big(\frac{t}{\sqrt{n}}\big) - e_m\big(\frac{t}{\sqrt{n}}\big)\bigg)\,  
e^{-t^2/2}\right] + r_n 
\en 
with 
\be 
|r_n| \leq \big(1 + |t|^{4m^2}\big)\, e^{-t^2/2}\  
\bigg(\frac{C}{n^{(m-1)/2}} + \frac{\ep_n}{n^{s-2}}\bigg). 
\en 
Here $C$, $c$ and $\ep_n$ are positive constants, depending on $s$ 
and the distribution of $X$, such that  
$\ep_n \rightarrow 0$, as $n \rightarrow \infty$. 
} 
 
\vskip5mm 
Thus, the closeness of $e_m$ to $v$ near zero determines 
the rate of approximation of $v_n$'s by the functions $u_m$'s 
(which have a different formal nature). This representation will be of  
use in the proof of Theorems 1.1-1.2, since, as we will see,  
the Liouville integrals may be applied to give a point wise bound  
on the (inverse)  
Fourier transforms within the class of the functions of the form  
$\hat V(\frac{t}{\sqrt{n}})\, e^{-t^2/2}$, such as in (1.9). 
 
Note that for $s \geq 3$ the expression in the last brackets of (1.10) is 
dominated by $C n^{-(m-1)/2}$, while in the range $2 \leq s < 3$ the 
second summand $\ep_n n^{-(s-2)}$ dominates the first one. In any case,  
with respect to the growing parameter $n$, the bound (1.9) is sharper  
than the one given in (1.8). This observation explains the improvement of  
(1.4) compared with the relation (1.3). 
 
We also remark that Theorem~1.3 holds for a more general class of functions,  
including Fourier-Stieltjes transforms $v(t)$ of finite (signed) measures  
with finite $s$-th moment, such that $v(0)=1$, $v'(0)=0$, $v''(0)=-1$. 
For example, the approximating functions $u_m$ and $e_m$ are not positive  
definite, but belong to this class. 
 
\vskip2mm 
The exposition of this paper is based on chapters of auxiliary results, which 
are organized in accordance with the following table. 
 
\vskip2mm 
{\it Contents} 
 
\vskip1mm 
1.\, Introduction. 
 
2.\, Differentiability with improved remainder terms. 
 
3.\, Differentiability of Fourier-Stieltjes transforms. 
 
4.\, Cumulants. The functions 
$\psi_z(t) = \frac{1}{2}\,t^2 + \frac{1}{z^2}\,\log v(tz)$. 
 
5.\, The case of moments of order $2 \leq s < 3$. 
 
6 \, Definition of the expansion polynomials $P_k$.  
 
7.\, Associated projection operators. 
 
8.\, Bounds of $P_k$ and their derivatives. 
 
9.\, Edgeworth-type expansion for the functions $v(tz)^{1/z^2}$. 
 
10. Proof of Theorem 1.3. 
 
11. Liouville fractional integrals and derivatives. 
 
12. Fourier transforms and fractional derivatives. 
 
13. Binomial decomposition of convolutions. 
 
14. Proof of Theorems 1.1 and 1.2.


\vskip5mm 
\section{{\bf Differentiability with Improved Remainder Terms}} 
\setcounter{equation}{0} 
 
\vskip2mm 
For our purposes, we use the following terminology. 
 
\vskip5mm 
{\bf Definition 2.1.} Let a complex-valued function $y = y(t)$ be  
defined in some interval $a<t<b$, and let $s \geq 0$.  
We say that $y$ is $s$-times differentiable, 
if it has continuous derivatives up to order $m = [s]$ in $(a,b)$, 
and for any $t_0 \in (a,b)$, as $t \rightarrow t_0$, 
\be 
y^{(m)}(t) = y^{(m)}(t_0) + o(|t-t_0|^{s-m}). 
\en 
 
\vskip2mm 
The case $s=0$ corresponds to continuity, while the case of a positive  
integer $s=m$ -- to the property of just having continuous derivatives  
up to order $m$. 
 
The following obvious characterization will be an important tool in the  
derivation of the Edgeworth-type expansions for charactersitic functions.  
It is obtained from (2.1) by the repeated integration over the variable  
$t$ near $t_0$.

\vskip5mm 
{\bf Proposition 2.2.} {\it Let $y$ have continuous derivatives of order up  
to $m = [s]$ in $(a,b)$. The function $y$ is $s$-times differentiable on  
$(a,b)$, if and only if, for any point $t_0 \in (a,b)$ and all $p = 0,\dots,m$,  
as $t \rightarrow t_0$, 
\be 
\frac{d^p}{dt^p}\, y(t) = \frac{d^p}{dt^p}\,  
\sum_{k=0}^{m} \frac{y^{(k)}(t_0)}{k!}\, (t-t_0)^k + o(|t-t_0|^{s-p}). 
\qquad 
\en 
} 
 
One can also provide quantitative estimates on the remainder term in (2.2),  
if we start with 
$$ 
|y^{(m)}(t) - y^{(m)}(t_0)| \leq |t-t_0|^{s-m} \ep(|t-t_0|), 
$$ 
where $\ep = \ep(w)$ is a non-decreasing function in $w \geq 0$, such 
that $\ep(w) \rightarrow 0$, as $w \rightarrow 0$. Then 
$$ 
\bigg|\frac{d^p}{dt^p}\, y(t) - \frac{d^p}{dt^p}\,  
\sum_{k=0}^{m} \frac{y^{(k)}(t_0)}{k!}\, (t-t_0)^k\bigg| \leq 
|t-t_0|^{s-p}\, \ep(|t-t_0|), 
$$ 
for any $p = 0,\dots,m$. 
 
\vskip2mm 
By the chain rule given below as Lemma 2.4, we have the following: 
 
\vskip5mm 
{\bf Proposition 2.3.} {\it If $y$ is $s$-times differentiable on $(a,b)$, $s \geq 0$, 
and $z = z(y)$ is analytic in some domain, containing all values $y(t)$, 
then the superposition $z(y(t))$ is also $s$-times differentiable on $(a,b)$. 
} 
 
\vskip5mm 
{\bf Lemma 2.4.} {\it Under the conditions of Proposition $2.3$, 
$z(y(t))$ has derivatives up to order $m = [s]$ on $(a,b)$, given by 
\be 
\frac{d^p}{dt^p}\ z(y(t)) = p! \sum  
\frac{d^{k_1 + \dots + k_p}\, z(y)}{dy^{k_1 + \dots + k_p}}\, \bigg|_{y=y(t)} 
\prod_{r=1}^p \frac{1}{k_r!}  
\bigg( \frac{1}{r!}\, \frac{d^r y(t)}{dt^r}\bigg)^{k_r}, 
\en 
for all $p = 1,\dots,m$, where the summation is performed over all non-negative integer 
solutions $(k_1,\dots,k_p)$ to the equation $k_1 + 2 k_2 + \dots + p k_p = p$. 
} 
 
\vskip5mm 
{\bf Proof of Proposition 2.3.} By the definition, for all  $t_0 \in (a,b)$ and 
$r = 1,\dots,m$, 
$$ 
\frac{d^r y(t)}{dt^r} = c_r + o(|t-t_0|^{s-m}), \quad  {\rm as} \ \ t \rightarrow t_0, 
$$ 
where $c_r = y^{(r)}(t_0)$. Raising these equalities  
to the $k_r$-powers and then multiplying them, we get a similar representation 
$$ 
\prod_{r=1}^m \frac{1}{k_r!}  
\bigg( \frac{1}{r!}\, \frac{d^r y(t)}{dt^r}\bigg)^{k_r} = c_k + o(|t-t_0|^{s-m}) 
$$ 
with some constants $c_k$, depending on the $m$-tuples $k = (k_1,\dots,k_m)$. 
In addition, putting $j = k_1 + \dots + k_m$, we have 
\bee 
z^{(j)}(y(t)) & = & z^{(j)}(y(t_0)) + O(|y(t) - y(t_0)|) \\ 
 & = & z^{(j)}(y(t_0)) + o(|t - t_0|^{s-m}). 
\ene 
Inserting these relations in (2.3) with $p=m$, we obtain that  
$[z(y(t))]^{(m)} = c + o(|t - t_0|^{s-m})$ with some constant $c$. 
 
In addition, the right-hand side of (2.3) represents a continuous function in $t$, 
so necessarily $c = [z(y)]^{(m)}(t_0)$. This means that (2.1) is fulfilled, 
and Proposition 2.3 is proved.

 
\vskip5mm 
\section{{\bf Differentiability of Fourier-Stieltjes Transforms}} 
\setcounter{equation}{0}

\vskip2mm 
A large variety of examples of $s$-times differentiable functions appear 
as Fourier-Stieltjes transforms of finite measures on the real line 
with finite absolute $s$-th moment.

\vskip5mm 
{\bf Proposition 3.1.} {\it Let $X$ be a random variable with characteristic  
function $v(t) = \E\, e^{itX}$. If $\E\, |X|^s < +\infty$, $s \geq 0$, then  
$v$ is $s$-times differentiable on the real line. Moreover, its $m = [s]$  
derivatives are representable, as $t \rightarrow 0$, by 
\be 
v^{(p)}(t) = \sum_{k=0}^{m-p} \E\, (iX)^{p+k}\, \frac{t^k}{k!} + o(|t|^{s-p}), 
\qquad p = 0,\dots,m. 
\en 
} 
 
One can state a similar proposition for more general Fourier-Stieltjes  
transforms 
$$ 
v(t) = \int_{-\infty}^{+\infty}\, e^{itx}\, dF(x), 
$$ 
where $F$ is a function of bounded variation on the real line, such that  
$\int |x|^s\,d\, |F|(x) < +\infty$ (where $|F|$ denotes  
the variation of $F$, viewed as a positive finite measure). 
On the other hand, such a more general statement may be obtained from  
Proposition~3.1, as well. Indeed, one can always represent $F$ as a linear  
combination $c_1 F_1 - c_2 F_2$ of two orthogonal probability distributions 
(with $c_1, c_2 \geq 0$). Then 
$|F| = c_1 F_1 + c_2 F_2$, so $\int |x|^s\,dF_i(x) < +\infty$. 
Applying Proposition 3.1 to $F_i$, we obtain a similar statement for $F$.

\vskip2mm 
{\bf Proof of Proposition 3.1.} By the moment assumption, the characteristic  
function $v$ has $m$ continuous derivatives, described by 
$$ 
v^{(p)}(t) = \E\, (iX)^p =  
\int_{-\infty}^{+\infty}\, e^{itx}\,(ix)^p\, dF(x), \qquad p = 0,1,\dots,m, 
$$ 
where $F$ is the distribution function of $X$. In particular, 
\be 
v^{(m)}(t) = \int_{-\infty}^{+\infty}\, e^{itx}\,(ix)^m\, dF(x). 
\en 
Hence, the relation (3.1) would follow immediately from  
the formula (2.2) of Proposition 2.2, once it has been established that $v$ is 
$s$-times differentiable. Namely, we need to see that, for any $t_0 \in \R$, as  
$t \rightarrow t_0$, 
\be 
v^{(m)}(t) = v^{(m)}(t_0) + o(|t-t_0|^{s-m}). 
\en 
 
The formula (3.2) is telling us that $v^{(m)}$ represents the Fourier-Stieltjes  
transform of the finite measure $F_m(dx) = (ix)^m\, dF(x)$ with  
$\int |x|^{s-m}\,d\, |F_m|(x) < +\infty$. Representing $F_m = c_1 G + c_2 H$  
with distribution functions $G$ and $H$, we can also represent $v^{(m)}$  
as a linear combination of the two Fourier-Stieltjes transforms with  
$\int |x|^{s-m}\,d\, G(x) < +\infty$ and similarly for $H$. 
Hence, in order to show that $v$ is $s$-times differentiable, 
it is enough to consider in (3.3) the case $m=0$ and $0 \leq s < 1$, only. 
 
Thus, Proposition 3.1 has been reduced to the case $\E\, |X|^s < +\infty$, $0 \leq s < 1$,  
when one needs to show that 
$$ 
v(t) = v(t_0) + o(|t-t_0|^s). 
$$ 
Moreover, without loss of generality, it suffices to consider the point $t_0 = 0$, 
only, in which case we need to show the relation $v(t) = 1 + o(|t|^s)$. 
The case $s=0$ is immediate, so let $s>0$ and write 
$$ 
1 - v(t) = \int_{-\infty}^{+\infty}\, (1-e^{itx})\,dF(x). 
$$ 
For definiteness, let $t>0$. Since in general $|1 - e^{ix}| \leq \min\{2,|x|\}$,  
$x \in \R$, we have 
\bee 
|1 - v(t)| & \leq &  
2 \int_{|x| \geq 2/t}\,dF(x) + t\int_{|x| < 2/t} |x|\,dF(x) \\ 
 & = &  
2 \int_{x \geq 2/t}\,dG(x) + t\int_{x < 2/t} x\,dG(x), \\ 
\ene 
where $G$ is the distribution of $|X|$. By the assumption, 
$$ 
\psi(x) = \E\, |X|^s\, 1_{\{|X| \geq x\}} = \int_x^{+\infty} y^s\,dG(y) 
\rightarrow 0, \quad {\rm as} \ \ x \rightarrow +\infty. 
$$ 
We have $\psi(x) \geq x^s \int_x^{+\infty} dG(y)$, so 
$\int_x^{+\infty} dG(y) = o(x^{-s})$, that is, 
$\int_{x \geq 2/t}\,dG(x) = o(t^s)$, as $t \rightarrow 0$. 
Finally, by integration by parts, 
$$ 
\int_{x < 2/t} x\,dG(x) \leq \int_0^{2/t} (1-G(x))\,dx 
\leq \int_0^{2/t} \frac{\psi(x)}{x^s}\,dx = 
t^{s-1} \int_0^2 \frac{\psi(y/t)}{y^s}\,dy. 
$$ 
Hence, 
$$ 
t \int_{x < 2/t} x\,dG(x) \leq t^s \int_0^2 \frac{\psi(y/t)}{y^s}\,dy. 
$$ 
But the last integral tends to zero, as $t \rightarrow 0$, by the 
Lebesgue dominated convergence theorem.  
 
Proposition 3.1 is proved.


\vskip10mm 
\section{{\bf Cumulants. The functions  
$\psi_z(t) = \frac{1}{2}\,t^2 + \frac{1}{z^2}\,\log v(tz)$}} 
\setcounter{equation}{0} 
 
\vskip2mm 
If a complex valued function $v$ on the real line has $m$ continuous derivatives 
and $v(0) \neq 0$, then $v(t) \neq 0$ in some interval $(-c,c)$. Moreover,  
the principal value of the logarithm $\log v(t)$ is well-defined in that  
interval and represents a function, which has also $m$ continuous derivatives.  
The corresponding derivatives at the origin (with a proper normalization), 
$$ 
\gamma_k = \frac{d^k}{i^k dt^k}\, \log v(t)\big|_{t=0}, \qquad k = 0,1,\dots,m, 
$$ 
will be called the generalized cumulants or just cumulants, associated to $v$. 
 
This terminology is standard, when $v$ represents the characteristic function 
of a random variable (with $m$ finite absolute moments). However, we shall have 
more general classes of functions. 
 
Applying Propositions 2.2-2.3, we arrive at: 
 
\vskip5mm 
{\bf Proposition 4.1.} {\it Let $v$ be $s$-times differentiable on the real line,  
$s \geq 0$, not vanishing in some interval $(-c,c)$. Then, $\log v$ is $s$-times  
differentiable in $(-c,c)$. In particular, as $t \rightarrow 0$, 
\be 
\frac{d^p}{dt^p}\, \log v(t) = \frac{d^p}{dt^p}\, 
\sum_{k=0}^{m} \frac{\gamma_k}{k!}\, (it)^k + o(|t|^{s-p}), \qquad 
p = 0,1,\dots,m, \ \ m = [s]. 
\en 
} 
 
\vskip2mm 
Note that if $v$ has $m+1$ continuous derivatives, then, by the usual Taylor's  
theorem, the remainder term in (4.1) can be sharpened, and we have 
\be 
\frac{d^p}{dt^p}\, \log v(t) = \frac{d^p}{dt^p}\, 
\sum_{k=0}^{m} \frac{\gamma_k}{k!}\, (it)^k + O(|t|^{(m+1)-p}), \qquad 
p = 0,1,\dots,m. 
\en 
 
If $v(t) = \E\, e^{itX}$ is the characteristic function of a random variable $X$, 
the assumptions of Proposition 4.1 are fulfilled, as long as 
$\E\, |X|^s < +\infty$ (Propositions 3.1). 
Then $\gamma_k$ are usual cumulants with $\gamma_0 = 0$, so (4.1) 
should be written 
$$ 
\frac{d^p}{dt^p}\, \log v(t) = \frac{d^p}{dt^p}\, 
\sum_{k=1}^{m} \frac{\gamma_k}{k!}\, (it)^k + o(|t|^{s-p}), \qquad p = 0,1,\dots,m. 
$$  
In the particular case $p = 0$, the above equality becomes 
\be 
\log v(t) = \sum_{k=1}^{m} \frac{\gamma_k}{k!}\, (it)^k + o(|t|^{s}). 
\en 
However, the case $p \geq 1$ in (4.1) cannot be deduced directly from (4.3). 
 
Let us recall how to relate the cumulants to the moments 
$\alpha_k = \E X^k$. Applying Lemma~2.4 with $z(y) = \log y$, $y = v(t)$, and 
at the point $t=0$, one obtains a well-known identity 
$$ 
\gamma_p\, =\, p!\, \sum\, (-1)^{k_1 + \dots + k_p - 1}\,  
(k_1 + \dots + k_p - 1)!\, 
\prod_{r=1}^p \frac{1}{k_r!}\, \bigg(\frac{\alpha_r}{r!}\bigg)^{k_r}, 
$$ 
where the summation extends over all non-negative integer solutions  
$(k_1,\dots,k_p)$ to the equation $k_1 + 2 k_2 + \dots + p k_p = p$. 
Note that $\gamma_p$ depends on the first $p$ moments of $X$, only. For example,  
\bee 
\gamma_1 & = & \alpha_1, \hskip23mm 
\gamma_3 \, = \, \alpha_3 - 3 \alpha_1 \alpha_2 + \alpha_1^3, \\ 
\gamma_2 & = & \alpha_2 - \alpha_1^2, \hskip14.5mm 
\gamma_4 \, = \, \alpha_4 - 3\alpha_3^2 - 4 \alpha_1 \alpha_3 +  
               12\, \alpha_1^2 \alpha_2 - 6 \alpha_1^4. 
\ene 
Under standard moment assumptions, such as $\E X = 0$, $\E X^2 = 1$,  
we have 
$\gamma_1 = 0$, $\gamma_2 = 1$, $\gamma_3 = \alpha_3$, $\gamma_4 = \alpha_4 - 3$. 
For any normal random variable, $\gamma_k = 0$, for all $k \geq 3$. 
 
\vskip5mm 
Now, returning to the general setting as in Proposition 4.1, assume that 
$s \geq 2$, and $v(0)=1$, $v'(0) = 0$, $v''(0)=-1$. Then  
$\gamma_0 = \gamma_1 = 0$, $\gamma_2 = 1$, so 
$\log v(t) = -\frac{t^2}{2} + o(|t|^2)$. 
Therefore, it is natural to center and normalize this function by introducing 
the family of the functions 
$$ 
\psi_z(t) = \frac{1}{2}\,t^2 + \frac{1}{z^2}\,\log v(tz), \qquad |tz| < c, 
$$ 
where $z \neq 0$ is a given parameter. Clearly, 
$\psi_z$ is $s$-times differentiable in this $t$-interval, and 
$$ 
\psi_z(0) = \psi_z'(0) = \psi_z''(0) = 0, \qquad  
\psi_z^{(p)}(0) = \gamma_p\ i^p z^{p-2} \ \ (p = 3,\dots,m). 
$$ 
Moreover, reformulating Proposition 4.1 in terms of the functions 
$t \rightarrow \frac{1}{z^2}\, \log v(tz)$ with fixed $z \neq 0$, we get:

\vskip5mm 
{\bf Corollary 4.2.} {\it Let $v(t)$ be $s$-times differentiable on the real line,  
$s \geq 2$, not vanishing for $|t| \leq c$, and such that $v(0)=1$, $v'(0) = 0$,  
$v''(0)=-1$. Given $z \neq 0$, in the interval $|tz| \leq c$, for all  
$p = 0,1,\dots,m$, $m = [s]$, 
\be 
\frac{d^p}{dt^p}\, \psi_z(t)  = \frac{d^p}{dt^p}\, 
\sum_{k=3}^{m} \frac{\gamma_{k}}{k!}\, (it)^{k}\,z^{k-2} +  
|t|^{s-p}\, |z|^{s-2} \ep(tz), 
\en 
where $\ep = \ep(t)$ is defined and continuous in $|t| \leq c$ and satisfies 
$\ep(t) \rightarrow 0$, as $t \rightarrow 0$. 
} 
 
\vskip5mm 
Also, as remarked after Proposition 4.1, cf. (4.2), when $v$ has $m+1$ continuous  
derivatives, a representation with sharper remainder term holds, 
$$ 
\frac{d^p}{dt^p}\, \psi_z(t)  = \frac{d^p}{dt^p}\, 
\sum_{k=3}^{m} \frac{\gamma_{k}}{k!}\, (it)^{k}\,z^{k-2} +  
A\, |t|^{(m+1)-p}\, |z|^{m-1}, 
$$ 
where $A = A(t,z)$ is a bounded function in the region $|tz| \leq c$. 
 
In the case of characteristic functions, Corollary 4.2 admits a slight refinement. 
 
\vskip5mm 
{\bf Corollary 4.3.} {\it Let $X$ be a random variable with the  
characteristic function $v$, and with $\E X = 0$, $\E X^2 = 1$,  
$\E\, |X|^s < +\infty$, $s \geq 2$. Then, given $z \neq 0$, the relation $(4.4)$ 
holds in the interval $|tz| < \sqrt{2}$. 
} 
 
\vskip5mm 
Indeed, by Taylor's theorem, $|1 - v(t)| \leq \frac{1}{2}\,t^2$, for all $t \in \R$.  
Hence, one may choose in Corollary 4.2 any value $0 < c < \sqrt{2}$.

 
\vskip10mm 
\section{{\bf The Case of Moments of Order $2 \leq s < 3$}} 
\setcounter{equation}{0} 
 
\vskip2mm 
In case $2 \leq s < 3$, Corollariy 4.2 is simplified, since then there 
are no terms in the sum (4.4). In particular, when $s=2$, we have  
\be 
\psi_z(t) = t^2\, \ep_0(tz), \quad \psi_z'(t) = |t|\, \ep_1(tz), \quad 
\psi_z''(t) = \ep_2(tz)  
\en 
with some functions $\ep_i(z) \rightarrow 0$, as $z \rightarrow 0$.  
This leads to the following observation which is classical in case  
of characteristic functions.

\vskip5mm 
{\bf Proposition 5.1.} {\it Assume that $v(t)$ has two continuous derivatives  
with $v(0)=1$, $v'(0) = 0$ and $v''(0)=-1$. There is a function  
$T_z \rightarrow +\infty$, for $z \rightarrow 0$ $(0 < |z| \leq 1)$, such that  
uniformly in the intervals $|t| \leq T_z$ 
\be 
\bigg|\frac{d^p}{dt^p}\, v(tz)^{1/z^2} - \frac{d^p}{dt^p}\, e^{-t^2/2}\bigg| 
 \, \leq \, e^{-t^2/2}\, \ep(z), \qquad p = 0,1,2, 
\en 
where $\ep(z) \rightarrow 0$, as $z \rightarrow 0$. 
} 
 
\vskip5mm 
{\bf Proof.} For completeness we include a well-known argument. 
Let $|t| \leq c$ be an interval where the function $v(t)$ is not vanishing. 
Choose the function $T = T_z$ as to satisfy $T_z |z| \leq c$, whenever  
$0<|z| \leq 1$, and $T_z |z| \rightarrow 0$, for $z \rightarrow 0$. 
These conditions will be assumed from now on. 
Moreover, for any continuous function $V(t)$, one 
can choose $T_z \rightarrow +\infty$, for $z \rightarrow 0$, such that 
$$ 
\sup_{|t| \leq T_z} |V(t) \psi_z(t)| \rightarrow 0, \quad {\rm as} \ \  
z \rightarrow 0, 
$$ 
and similarly for the first two derivatives of $\psi_z$. 
 
For the proof, it is enough to see that, whenever $\ep(z) \rightarrow 0$ 
and $W(t) \geq 0$ is continuous and increasing in $t \geq 0$, one can choose  
$T_z \rightarrow +\infty$, such that 
$W(T_z)\, \sup_{|t| \leq T_z} |\ep(tz)| \rightarrow 0$, as $z \rightarrow 0$. 
Here, we may assume in the following that $\ep(z) \geq 0$ is even and also  
increasing in $z > 0$. Then, the latter statement may be simplified to  
$W(T_z)\,\ep(T_z z) \rightarrow 0$, which is obviously true with a sufficiently 
slowly growing $T_z$. 
 
In particular, in view of (5.1), with some $T_z \rightarrow +\infty$, as  
$z \rightarrow 0$ $(0 < |z| \leq 1)$, we have 
\be 
\ep(z) = 
\sup_{|t| \leq T_z} \left(|\psi_z(t)| + |\psi_z'(t)| + |\psi_z''(t)|\right)  
\rightarrow 0, \quad {\rm as} \ \ z \rightarrow 0. 
\en 
 
Now, write $v(tz)^{1/z^2} = g(t) e^{\psi_z(t)}$, where 
$g(t) = e^{-t^2/2}$. Applying (5.3), we get 
$$ 
|v(tz)^{1/z^2} - g(t)| \leq g(t)\,|e^{\psi_z(t)} - 1| \leq Cg(t)|\psi_z(t)|, 
$$ 
for all $|t| \leq T_z$ with some constant $C$. Since also  
$\psi_z(t) \rightarrow 0$ uniformly in that interval,  
we arrive at the desired conclusion in case $p=0$.  
 
Writing 
$(v(tz)^{1/z^2})' = g'(t) e^{\psi_z(t)} + g(t) \psi_z'(t) e^{\psi_z(t)}$ 
and using the previous step, we get 
$$ 
|(v(tz)^{1/z^2})' - g'(t)| \leq C|t| g(t)\,|\psi_z(t)| + C g(t)\,|\psi_z'(t)|. 
$$ 
Since $\psi_z'(t) \rightarrow 0$ and $t\psi_z(t) \rightarrow 0$ uniformly  
in that interval with an appropriate choice of $T_z$, we arrive at the 
conclusion in case $p=1$. The case $p=2$ is similar. 
 
\vskip5mm 
Now, let us turn to the range $2 < s < 3$. In this case, we obtain up to  
polynomial factors in front of $e^{-t^/2}$ in (5.2) more 
information about the possible growth of $T_z$. 
 
\vskip5mm 
{\bf Proposition 5.2.} {\it Let $v(t)$ be $s$-times differentiable, $2 < s < 3$,  
not vanishing for $|t| \leq c$ $(c > 0)$, and such that $v(0)=1$, $v'(0) = 0$,  
$v''(0)=-1$. Given $0 < |z| \leq 1$,  
$$ 
\bigg|\frac{d^p}{dt^p}\, v(tz)^{1/z^2} - \frac{d^p}{dt^p}\, e^{-t^2/2}\bigg| 
 \, \leq \, \left(|t|^{s-2} + |t|^{s+2}\right) e^{-t^2/2}\, \ep(z), \qquad p = 0,1,2, 
$$ 
uniformly for $|t| \leq c\,|z|^{-(s-2)/s}$ 
with some function $\ep(z) \rightarrow 0$, as $z \rightarrow 0$. 
} 
 
\vskip5mm 
{\bf Proof.} By Corollary 4.2, in the intervals  
$|t| \leq T_z = c\,|z|^{-(s-2)/s}$ with $0 < |z| \leq 1$ 
\bee 
|\psi_z(t)| & \leq & |t|^s\ |z|^{s-2} \ \ep(z), \\ 
|\psi_z'(t)| & \leq & |t|^{s-1}\, |z|^{s-2}\, \ep(z), \\ 
|\psi_z''(t)| & \leq & |t|^{s-2}\, |z|^{s-2}\, \ep(z) 
\ene 
with some bounded function $\ep(z) \rightarrow 0$, as $z \rightarrow 0$. 
Indeed, the conditions $0 < |z| \leq 1$ and $|t| \leq T_z$ 
insure that $|tz| \leq c\,|z|^{2/s} \leq c$ and also $|tz| \rightarrow 0$, as 
$z \rightarrow 0$, uniformly in $|t| \leq T_z$. 
 
Now, by the first inequality, 
$$ 
|\psi_z(t)| \leq |t|^s\ |z|^{s-2} \, \ep(z) \leq C \qquad (|t| \leq T_z) 
$$ 
with some constant $C$. Hence, using the same notations and arguments as in the  
proof of Proposition 5.1, for all $|t| \leq T_z$, 
$$ 
|v(tz)^{1/z^2} - g(t)| \leq C' g(t)|\psi_z(t)| \leq C' g(t)  
\cdot |t|^s\ |z|^{s-2} \, \ep(z). 
$$ 
Since also $|\psi_z'(t)| \leq |t|^{s-1}\ |z|^{s-2} \, \ep(z)$, with some  
constants $C,C'$ we get 
\bee 
|(v(tz)^{1/z^2})' - g'(t)| 
 & \leq &  
C|t| g(t)\,|\psi_z(t)| + C g(t)\,|\psi_z'(t)| \\ 
 & \leq & 
C' \left(|t|^{s+1} + |t|^{s-1}\right)\, g(t)\,\ep(z). 
\ene 
Finally, writing 
$$ 
(v(tz)^{1/z^2})'' = g''(t) e^{\psi_z(t)} + 2 g'(t) \psi_z'(t) e^{\psi_z(t)} 
+ 2 g(t) (\psi_z'' + \psi_z'(t)^2)\, e^{\psi_z(t)} 
$$ 
and using $|\psi_z''(t)| \leq |t|^{s-2}\ |z|^{s-2} \, \ep(z)$, 
we get that up to some constants 
\bee 
\left|(v(tz)^{1/z^2})'' - g''(t)\right| 
 & \leq &  
Ct^2 g(t)\,|\psi_z(t)| + 
C|t| g(t)\,|\psi_z'(t)| + C g(t)\,\left(|\psi_z''(t)| + |\psi_z'(t)|^2\right) \\ 
 & \leq & 
C' \left(|t|^{s+2} + |t|^s + |t|^{s-2} + |t|^{2(s-1)}\right)\, g(t)\  
|z|^{s-2}\ep(z). 
\ene 
All powers of $|t|$ vary from $s-2$ to $s+2$, so Proposition 5.2 is proved.

 
\vskip10mm 
\section{{\bf Definition of the Expansion Polynomials $P_k$}} 
\setcounter{equation}{0} 
 
\vskip2mm 
The polynomials $P_k$ introduced in Section 1 appear not only in connection  
with charactersitic functions, but in a more general setting as well.  
 
Namely, let $v(t)$ be a complex-valued function on the real line, which is 
$s$-times differentiable $(s \geq 2)$ and such that 
$v(0)=1$, $v'(0)=0$, $v''(0)=-1$.  
Then $v$ has cumulants $\gamma_k$, $k = 1,\dots,m$, where $m = [s]$. 
Moreover, $\gamma_1 = 0$ and $\gamma_2 = 1$. 
 
Assume $v$ is not vanishing in the interval $|t| \leq c$, and let us return  
to the functions  
$$ 
\psi_z(t) = \frac{1}{2}\,t^2 + \frac{1}{z^2}\,\log v(tz), 
$$ 
where $z \neq 0$ is viewed as a (small) parameter and $|tz| \leq c$. 
Recall that, by Corollary 4.2, 
\be 
e^{t^2/2}\, v(tz)^{1/z^2} = e^{\psi_z(t)} = \exp 
\bigg\{\sum_{k=1}^{m-2} \frac{\gamma_{k+2}}{(k+2)!}\, (it)^{k+2}\,z^k +  
|t|^{s} |z|^{s-2}\, \ep(tz)\bigg\}, 
\en 
where $\ep(t)$ is defined and continuous in $|t|\leq c$ and satisfies 
$\ep(t) \rightarrow 0$, as $t \rightarrow 0$. Moreover, if $v$ has $m+1$  
derivatives, the remainder term here may be replaced with 
$A\, |t|^{(m+1)-p}\, |z|^{m-1}$, where $A =A(t,z)$ is bounded in the region  
$|tz| \leq c$. 
 
The sum in (6.1) is vanishing in case $2 \leq s < 3$. To study 
this representation in the case $s \geq 3$, introduce the polynomials 
$$ 
W_z(t) = \sum_{k=1}^{m-2} \frac{\gamma_{k+2}}{(k+2)!}\, (it)^{k+2}\,z^k. 
$$ 
By a formal Taylor's representation with respect to the (complex) variable $z$,  
we have 
\be 
e^{W_z(t)} = 1 + \sum_{k=1}^\infty a_k(it)\,z^k 
\en 
with some coefficients $a_k(it)$. 
To justify this step and precisely determine the coefficients, write 
$$ 
e^{W_z(t)} \, = \, \sum_{p=0}^\infty \frac{W_z(t)^p}{p!}, 
$$ 
which can further be expanded as 
\be 
\sum_{p=0}^\infty \ \ \sum_{p_1 + \dots + p_{m-2} = p}   
\frac{(\frac{\gamma_3}{3!})^{p_1} \dots 
(\frac{\gamma_{m}}{m!})^{p_{m-2}}}{p_1!\dots p_{m-2}!}\, 
(it)^{3p_1 + 4p_2 + \dots + m p_{m-2}}\, z^{p_1 + 2 p_2 + \dots + (m-2) p_{m-2}}. 
\en 
The whole double sum is absolutely summable for all complex numbers 
$t$ and $z$. Indeed, let $C = \sum_{k = 3}^m \frac{|\gamma_k|}{k!}$. 
For any fixed integer $p \geq 0$, the finite sum of the absolute values 
of the terms  in (6.3) is bounded by 
$$ 
\frac{1}{p!}\, 
\bigg(\sum_{k=1}^{m-2} \frac{|\gamma_{k+2}|}{(k+2)!}\, |t|^{k+2}\,|z|^k\bigg)^p 
\leq \frac{1}{p!}\,C^p \left(\max_{1 \leq k \leq m-2} |t|^{k+2}\,|z|^k\right)^p. 
$$ 
Assume without loss of generality (in order to get some quantitative bounds) that 
$|t^3 z| \leq 1$ and $|z| \leq 1$. Then, 
$$ 
|t|^{k+2}\,|z|^k \leq |z|^{-(k+2)/3}\,|z|^k = |z|^{(2k-2)/3} \leq 1. 
$$ 
Hence, the above sum is bounded by $C^p/p!$\, 
Furthermore, note that $|W_z(t)| \leq C$.

Thus, the total sum of the absolute values is bounded by $e^C$, and 
one may freely choose the order of summation. Collecting the coefficients  
in (6.3) in front of $z^k$, we arrive at (6.2) with 
$$ 
a_k(it) \ \ = \sum_{p_1 + 2 p_2 + \dots + (m-2) p_{m-2} = k}   
\frac{1}{p_1!\dots p_{m-2}!}\, 
\bigg(\frac{\gamma_3}{3!}\bigg)^{p_1} \dots 
\bigg(\frac{\gamma_{m}}{m!}\bigg)^{p_{m-2}} 
(it)^{3p_1 + 4p_2 + \dots + m p_{m-2}}, 
$$ 
where the summation is extended over all non-negative integer  
solutions $(p_1,\dots,p_{m-2})$ to the equation  
$p_1 + 2 p_2 + \dots + (m-2) p_{m-2} = k$. Note that 
$$ 
3p_1 + 4p_2 + \dots + m p_{m-2} \, = \, k + 2(p_1 + p_2 + \dots + p_{m-2}). 
$$ 
Hence, replacing $it$ with $t$, 
\be 
a_k(t) \ \ = \sum_{p_1 + 2 p_2 + \dots + (m-2) p_{m-2} = k}   
\frac{1}{p_1!\dots p_{m-2}!}\, 
\bigg(\frac{\gamma_3}{3!}\bigg)^{p_1} \dots 
\bigg(\frac{\gamma_{m}}{m!}\bigg)^{p_{m-2}} 
t^{k + 2(p_1 + \dots + p_{m-2})}. 
\en 
 
In addition, if $k \leq m-2$, the condition $p_1 + 2 p_2 + \dots + (m-2) p_{m-2} = k$ 
implies that $p_{k+1} = \dots = p_{m-2} = 0$. Therefore, in this case $a_k$ depends  
on the first $k$ cumulants $\gamma_3,\dots,\gamma_{k+2}$, only. More precisely, 
$$ 
a_k(t) \ = \sum_{p_1 + 2 p_2 + \dots + k p_k = k}   
\frac{1}{p_1!\dots p_k!}\, 
\bigg(\frac{\gamma_3}{3!}\bigg)^{p_1} \dots 
\bigg(\frac{\gamma_{k+2}}{(k+2)!}\bigg)^{p_k} 
t^{k + 2(p_1 + \dots + p_k)}, \ \ 1 \leq k \leq m-2. 
$$ 
 
For example, if $m=3$, we have $a_1(it) = \frac{\gamma_3}{6}\,(it)^3$. 
In case $m=4$, $a_1(t) = \frac{\gamma_3}{6}\,t^3$ and 
$$ 
a_2(t) = \frac{\gamma_3^2}{72}\,t^6 + \frac{\gamma_4}{24}\,t^4. 
$$ 
 
In general, subject to $p_1 + 2 p_2 + \dots + k p_k = k$, the expression  
$k + 2(p_1 + \dots + p_k)$ does not exceed $3k$ and reaches this value 
(when $p_1 = k$, $p_2 = \dots = p_k = 0$), 
so $a_k$ represents a polynomial in $t$ of degree exactly $3k$. 
 
\vskip5mm 
{\bf Definition 6.1.} Given an integer $m \geq 3$ and complex numbers 
$\gamma_3,\dots,\gamma_m$, one defines $P_k$ $(1 \leq k \leq m-2)$ 
as the polynomial $a_k$ introduced above, namely, 
$$ 
P_k(t) \ \ = \sum_{p_1 + 2 p_2 + \dots + k p_k = k}   
\frac{1}{p_1!\dots p_k!}\, 
\bigg(\frac{\gamma_3}{3!}\bigg)^{p_1} \dots 
\bigg(\frac{\gamma_{k+2}}{(k+2)!}\bigg)^{p_k} 
t^{k + 2(p_1 + \dots + p_k)}. 
$$ 
 
\vskip2mm 
With this definition, the representation (6.2) may also be written as 
\be 
\exp\bigg\{\sum_{k=3}^{m} \frac{\gamma_k}{k!}\,(it)^k z^{k-2}\bigg\}\, = \, 
1 + \sum_{k=1}^{m-2} P_k(it)\,z^k + \sum_{k=m-1}^\infty a_k(it) z^k, 
\en 
where $a_k$'s are described in (6.4).

 
\vskip10mm 
\section{{\bf Associated Projection Operators}} 
\setcounter{equation}{0} 
 
\vskip2mm 
Let us note first that every polynomial $a_k$ in (6.4) contains terms involving  
powers of $t$ not smaller than $k+2$ (since necessarily  
$p_1 + \dots + p_{m-2} \geq 1$). This observation may be used to obtain  
an initial trivial bound for the last sum in (6.5) in case of small values of $t$. 
 
\vskip5mm 
{\bf Lemma 7.1.} {\it Given complex numbers $\gamma_3,\dots,\gamma_m$, for all 
complex $z$ and $t$, $|t| \leq 1$, 
\be 
\sum_{k=m-1}^\infty |a_k(it) z^k| \leq \big(e^{C(z)} - 1\big)\, |t|^{m+1}, 
\en 
where $C(z) = \sum_{k=3}^{m} \frac{|\gamma_k|}{k!}\,|z|^{k-2}$. 
} 
 
\vskip5mm 
Indeed, using (6.4) and $|t| \leq 1$, we see that each term in the sum of (7.1)  
is bounded by $|t|^{m+1}$ up to the factor 
$$ 
\sum_{p_1 + 2 p_2 + \dots + (m-2) p_{m-2} = k} \frac{1}{p_1!\dots p_{m-2}!}\, 
\bigg(\frac{|\gamma_3|\, |z|}{3!}\bigg)^{p_1} \dots 
\bigg(\frac{|\gamma_{m}|\, |z|^{m-2}}{m!}\bigg)^{p_{m-2}}. 
$$ 
Summation of these expressions over all $k \geq 1$ results in $e^{C(z)} - 1$. 
 
The bound of Lemma 7.1 is needed in order to express the "cumulants" $\gamma_k$ 
directly in terms of the associated polynomials $P_j$.

\vskip5mm 
{\bf Lemma 7.2.} {\it Given complex numbers $\gamma_3,\dots,\gamma_m$, for all 
complex $z$ and $k=3,\dots,m$, 
\be 
\frac{\gamma_k}{k!}\, z^{k-2}\, =\, \frac{d^k}{i^k\, dt^k}\, 
\log\bigg(1 + \sum_{j=1}^{m-2} P_j(it) z^j\bigg)\bigg|_{t=0}. 
\en 
} 
 
\vskip2mm 
To see this, note that, by (6.5) and (7.1), we obtain that, as $t \rightarrow 0$, 
$$ 
\sum_{k=3}^{m} \frac{\gamma_k}{k!}\,(it)^k z^{k-2} =  
\log\bigg(1 + \sum_{j=1}^{m-2} P_j(it) z^j\bigg) + O(|t|^{m+1}). 
$$ 
Since the right-hand side (without the remainder term) represents an analytic  
function near zero, comparison of coefficients powers of $t$ immediately leads  
to (7.2). 
 
The identity (6.5) suggests to consider special operators, defined on the space 
$V_s$ of all $s$-times differentiable functions $v:\R \rightarrow \C$, $s \geq 2$,  
such that $v(0)=1$, $v'(0)=0$, $v''(0)=-1$. 
 
\vskip5mm 
{\bf Definition 7.3.} Given $v \in V_s$, $s \geq 2$, and an integer $2 \leq m \leq s$,   
we put 
$$ 
(T_m v)(t) = e^{-t^2/2} \bigg(1 + \sum_{k=1}^{m-2} P_k(it)\bigg), \quad t \in \R, 
$$ 
where $P_k$ are the polynomials from Definition 6.1 based on the cumulants 
$$ 
\gamma_k = \frac{d^k}{i^k\, dt^k}\,\log v(t)\big|_{t=0}, \qquad k = 3,\dots,m-2. 
$$ 
 
\vskip2mm 
If $m=2$, there are no cumulants and polynomials in the definition, so 
$(T_2 v)(t) = e^{-t^2/2}$, for any $v \in V_s$. If $m=3$, 
$$ 
(T_3 v)(t) = e^{-t^2/2} \bigg(1 + \frac{\gamma_3}{6}\,(it)^3\bigg), 
$$ 
for any $v \in V_s$, $s \geq 3$ (where $\gamma_3$ may be an arbitrary complex number). 
If $m=4$, then for any $v \in V_s$, $s \geq 4$, 
$$ 
(T_4 v)(t) = e^{-t^2/2} \bigg(1 + \frac{\gamma_3}{6}\,(it)^3 + 
\frac{\gamma_4}{24}\,(it)^4 + \frac{\gamma_3^2}{72}\,(it)^6\bigg). 
$$

Clearly, every $T_m v$ is an entire function and hence belongs to all $V_s$, 
$s \geq m$. This defines an operator $T_m: V_s \rightarrow V_s$ which turns out 
a projection operator

\vskip5mm 
{\bf Proposition 7.4.} {\it We have $T_m T_m v = T_m v$ for any $v \in V_s$, 
$2 \leq m \leq s$. Moreover, $T_m v$ and $v$ have identical derivatives  
at the origin up to order $m$. 
} 
 
\vskip5mm 
{\bf Proof.} The statement is equivalent to the property that $e_m = T_m v$ and $v$  
have equal cumulants. Let $\widetilde\gamma_k$ and $\gamma_k$ denote  
the cumulants of $e_m$ and $v$, respectively ($3 \leq k \leq m$). By Definition 7.3, 
$$ 
\frac{\widetilde\gamma_k}{k!}\, =\,  
\frac{d^k}{i^k\, dt^k}\,\log e_m(t)\big|_{t=0} = 
\frac{d^k}{i^k\, dt^k}\, 
\log\bigg(1 + \sum_{j=1}^{m-2} P_j(it)\bigg)\bigg|_{t=0}. 
$$ 
But the right-hand side equals $\frac{\gamma_k}{k!}$, according to Lemma 7.2 
applied with $z = 1$. 
 
Thus Proposition 7.4 is proved. 
 
\vskip5mm 
Note that $T_m v$ does not need to be a characteristic function, even if 
$v$ is a characteristic function of some random variable. However, it always 
represents the Fourier-Stieltjes transform of a finite signed measure. 
 
In the following we approximate $v$ by its projections $T_m v$. 
Combining Proposition 7.4 with Proposition 2.2, we get: 
 
\vskip5mm 
{\bf Corollary 7.5.} {\it Given $v \in V_s$, $2 \leq m \leq s$, as $t \rightarrow 0$, 
$$ 
\frac{d^p}{dt^p}\, \big(v(t) - T_m v(t)\big) = o(|t|^{s-p}), \qquad 
p = 0,1,\dots,m. 
$$ 
}

Finally, let us formulate an asymptotic property of the projection operators $T_m$ 
for growing parameter $m$ (although this will not be needed in the sequel).

\vskip5mm 
{\bf Proposition 7.6.} {\it Assume that $v(t)$ admits an analytic extension to the disc 
$|t| < \rho$, where it has no zeros, and $v(0)=1$, $v'(0)=0$, $v''(0)=-1$. Then 
$T_m v(t) \rightarrow v(t)$, as $m \rightarrow \infty$, i.e., 
$$ 
v(t) = e^{-t^2/2} \bigg(1 + \sum_{k=1}^{\infty} P_k(it)\bigg), \qquad |t| < \rho. 
$$ 
Moreover, the series is convergent absolutely. 
 
} 
 
\vskip5mm 
If $v(t) = \E\, e^{itX}$ is the characteristic function of a random variable $X$, 
the assumptions of Proposition 7.6 are fulfilled, provided that 
$\E X = 0$, $\E\, X^2 = 1$, $\E\, e^{\rho |X|} < +\infty$ 
(that is, an exponential moment of order $\rho$ is finite) and $v(t)$ 
does not vanish in the disc $|t|< \rho$. 
 
\vskip5mm 
{\bf Proof.} By the assumption, $\log v(t)$ is analytic in the disc $|t| < \rho$,  
so it is representable as the sum of the absolutely convergent power series 
\be 
\log v(t) = \sum_{k=3}^{\infty} \frac{\gamma_k}{k!}\,(it)^k, \qquad |t| < \rho. 
\en 
Hence, starting with (6.5) with $z=1$ and letting there $m \rightarrow \infty$,  
it is sufficient to show that $\sum_{k=m-1}^\infty |a_k(it)| \rightarrow 0$ 
(note that $a_k$'s also depend on $m$). 
 
Rewrite the representation (6.4) as 
$$ 
a_k(t) \ = \sum_{p_1 + 2 p_2 + \dots + (m-2) p_{m-2} = k}   
\frac{1}{p_1!\dots p_{m-2}!}\, \bigg(\frac{\gamma_3\, t^3}{3!}\bigg)^{p_1} \dots 
\bigg(\frac{\gamma_m\, t^m}{m!}\bigg)^{p_{m-2}}, 
$$ 
which implies that 
$$ 
|a_k(t)| \ \leq \sum_{p_1 + 2 p_2 + \dots + (m-2) p_{m-2} = k}   
\frac{1}{p_1!\dots p_{m-2}!}\, \bigg(\frac{|\gamma_3|\, |t|^3}{3!}\bigg)^{p_1} \dots 
\bigg(\frac{|\gamma_m|\, |t|^m}{m!}\bigg)^{p_{m-2}}. 
$$ 
Here the right-hand side may be bounded by the quantity 
$$ 
b_k(t) \ = \sum_{p_1 + 2 p_2 + 3 p_3 + \dots = k} \, \prod_{r=1}^\infty 
\frac{1}{p_r!}\, \bigg(\frac{|\gamma_{r+2}|\, |t|^{r+2}}{(r+2)!}\bigg)^{p_r}, 
$$ 
which does not depend on $m$. 
After summation over all $k \geq 1$ (thus removing any constraint on $p_r$),  
we get $\sum_{k=1}^\infty b_k(t) = e^{C(|t|)} - 1$, where 
$C(a) = \sum_{k=3}^{\infty} \frac{|\gamma_k|}{k!}\,a^k$. 
But $C(|t|) < +\infty$ for all $|t| < \rho$ in view of the absolute convergence 
of the series (7.3). Hence, in this case 
$$ 
\sum_{k=m-1}^\infty |a_k(it)| \leq \sum_{k=m-1}^\infty b_k(t) \rightarrow 0, \quad 
{\rm as} \ \ m \rightarrow \infty. 
$$ 
 
With similar arguments, we also obtain that $\sum_{k=1}^\infty |P_k(it)| < +\infty$  
for $|t| < \rho$, in view of Definition 6.1. Thus Proposition 7.6 is proved.


\vskip10mm 
\section{{\bf Bounds of $P_k$ and their Derivatives}} 
\setcounter{equation}{0} 
 
\vskip2mm 
We will need a bound similar to the one in Lemma 7.1 which extends to large 
values of $t$ and involving derivatives of the polynomials $a_k$ and $P_k$. 
 
To this aim, we start with arbitrary complex numbers 
$\gamma_3,\dots,\gamma_m$, $m \geq 3$ (which may be interpreted as 
cumulants of a given function $v$) and return to the representation (6.2), 
\be 
w_z(t) = e^{W_z(t)} = 1 + \sum_{k=1}^\infty a_k(it)\,z^k, \qquad t,z \in \C, 
\en 
where 
\be 
W_z(t) = \sum_{k=1}^{m-2} \frac{\gamma_{k+2}}{(k+2)!}\, (it)^{k+2}\,z^k, 
\en 
and where the polynomials $a_k$ are described in (6.4). 
By the very definition, $a_k = P_k$, as long as $k \leq m-2$. 
 
As we have already noticed, the sum in (8.1) is absolutely convergent  
and therefore represents an entire function with respect to $z$,  
for any fixed $t$.  
It is also clear that the series may be term wise differentiated, so that 
\be 
w_z^{(p)}(t) = \sum_{k=1}^\infty i^p\, a_k^{(p)}(it)\,z^k, \qquad p \geq 1, 
\en 
which is absolutely convergent, as well.  
 
In order to bound $a_k$ and its derivatives, we use the quantity 
$$ 
C = \sum_{k = 3}^m |\gamma_k|. 
$$ 
One natural approach (which is however different than the one in [Pe2]) 
is based on the application of Cauchy's integral formula 
$$ 
a_k(it) = \frac{1}{2\pi i} \int_{|z| = \rho} \frac{w_z(t)}{z^{k+1}}\,dz 
$$ 
with a suitably chosen parameter $\rho > 0$. In view of (8.3), there is  
a more general identity, involving the derivatives, 
\be 
i^p\, a_k^{(p)}(it) = \frac{1}{2\pi i} \int_{|z| = \rho}  
\frac{w_z^{(p)}(t)}{z^{k+1}}\,dz, \qquad p = 0,1,2\dots 
\en

\vskip5mm 
{\bf Lemma 8.1.} {\it Let $|tz| \leq 2$ and $|t^3 z| \leq 2$, and let 
$0 \leq p \leq m$ be an integer. Then 
\be 
|W_z^{(p)}(t)| \leq 2^{m-2} C\, |t|^{-p} \min\{1,|t|^2\}. 
\en 
} 
 
\vskip2mm 
Indeed, by definition (8.2),  
$$ 
W_z^{(p)}(t) =  
\sum_{q = \max(p,3)}^{m} \frac{\gamma_q\, i^q}{(q-p)!}\, t^{q-p} z^{q-2}. 
$$ 
But 
$|t^{q-p}\,z^{q-2}| = |tz|^{q-3}\, |t^3 z|\, |t|^{-p} \leq 2^{m-2}\, |t|^{-p}$,  
whenever $3 \leq q \leq m$. Hence, $|W_z^{(p)}(t)| \leq 2^{m-2} C \, |t|^{-p}$. 
On the other hand, just using $|z| \leq 2/|t|$, we get 
$$ 
|t^{q-p}\,z^{q-2}| \leq |t|^{q-p}\, \frac{2^{q-2}}{|t|^{q-2}} =  
2^{q-2} |t|^{2 - p} \leq 2^{m-2}  |t|^{2 - p}. 
$$ 
This gives an improvement over the previous estimate in case $|t| \leq 1$ 
and proves (8.5).

\vskip5mm 
{\bf Lemma 8.2.} {\it For all integers $k \geq 1$, $0 \leq p \leq m$, and 
all complex $t$, 
$$ 
|a_k^{(p)}(it)| \leq C_{m,p}\,|t|^{-p} \min\{1,|t|^2\}\, 
\bigg(\frac{\max\{|t|,|t|^3\}}{2}\bigg)^k 
$$ 
with constants $C_{m,p} = (4^m (1+C))^p\ e^{2^m C}$. 
} 
 
\vskip5mm 
{\bf Proof.} Given $t \neq 0$, we choose in (8.4) the radius 
$$ 
\rho = \frac{2}{\max\{|t|,|t|^3\}}.  
$$ 
Hence, on the circle $|z|=\rho$, both $|tz| \leq 2$ and $|t^3 z| \leq 2$  
are fulfilled, thus the inequality (8.5) may be applied. In particular,  
$|W_z(t)| \leq 2^{m-2} C$, and from (8.4) with $p=0$ we get the desired  
estimate 
$$ 
|a_k(it)| \leq \frac{1}{\rho^k}\, e^{2^{m-2} C}. 
$$ 
 
Next, by the formula (2.3) of Lemma 2.4, for all $p \geq 1$, 
$$ 
w_z^{(p)}(t) = p!\, w_z(t) \sum \prod_{r=1}^p \frac{1}{k_r!}  
\bigg(\frac{W_z^{(r)}(t)}{r!}\bigg)^{k_r}, 
$$ 
where the summation is taken over all non-negative integer solutions  
$(k_1,\dots,k_p)$ to the equation $k_1 + 2 k_2 + \dots + p k_p = p$.  
Hence, using (8.5), given that $|z| = \rho$, we arrive at 
$$ 
|w_z^{(p)}(t)| \, \leq \, e^{2^{m-2} C} |t|^{-p}\, p!\, \sum \prod_{r=1}^p \frac{1}{k_r!}  
\bigg(\frac{2^{m-2} C\,\min\{1,|t|^2\}}{r!}\bigg)^{k_r}. 
$$ 
Since necessarily $1 \leq k_1 + \dots + p k_p \leq p$, the product may be 
bounded by the product of $\min\{1,|t|^2\}$ (in the first power) 
and $2^{m-2} C$ (replaced by $2^{m-2} (1+C)$), and raised to 
power $p$. This leads to 
\be 
|w_z^{(p)}(t)| \, \leq \, e^{2^{m-2} C} \left(2^{m-2} (1+C)\right)^p  
|t|^{-p} \min\{1,|t|^2\}\, B_p, 
\en 
where 
$B_p = p! \sum \prod_{r=1}^p \frac{1}{k_r!} (\frac{1}{r!})^{k_r}$. 
This constant can also be described by virtue of the same formula (2.3),  
applied with $z = e^y$ and $y(s) = e^{s}$, in which case it reads 
$$ 
\frac{d^p}{ds^p}\ e^{e^s} = e^{e^s}\, p! \sum  
\prod_{r=1}^p \frac{1}{k_r!} \bigg(\frac{e^s}{r!}\bigg)^{k_r}. 
$$ 
One should apply this formula at $s = 0$, thus we consider the functions 
$b_p(s)  = (e^{e^s})^{(p)}$, $s \geq 0$, and their values  
$b_p = b_p(0) = B_p/e$. The recursive identity  
$b_{p+1}(s) = (e^s\,e^{e^s})^{(p)} = e^s \sum_{r=0}^p C_p^r\, b_r(s)$  
implies that the sequence $r \rightarrow b_r$ is non-decreasing and 
$b_{p+1} \leq 2^p\, b_p$. Therefore, 
$$ 
b_p \leq 2^{p-1} b_{p-1} \leq 2^{p-1} 2^{p-2} b_{p-2} \leq \dots \leq 
2^{p-1} 2^{p-2} \dots\, 2^0\, b_0 = 2^{p(p-1)/2}\,e. 
$$ 
Hence, $B_p \leq 2^{p(p-1)/2}$, and together with (8.6) this gives the estimate  
$$ 
|w_z^{(p)}(t)| \, \leq \, e^{2^{m-2} C} \big(2^{(m-2) + (p-1)/2} (1+C)\big)^p  
|t|^{-p} \min\{1,|t|^2\}. 
$$ 
It remains to apply (8.4) and simplify the constant.  
Thus Lemma 8.2 is proved.

\vskip5mm 
Now, fix an integer $p = 0,1\dots,m$, and assume that  
$|z| \leq \frac{\rho}{2} = \frac{1}{\max\{|t|,|t|^3\}}$, that is,  
$|tz| \leq 1$ and $|t^3 z| \leq 1$. By Lemma~8.2, 
\bee 
\sum_{k=m-1}^\infty |a_k^{(p)}(it)\,z^k| 
 & \leq &  
C_{m,p}\,|t|^{-p} \min\{1,|t|^2\}\, \sum_{k=m-1}^\infty \frac{|z|^k}{\rho^k} \\ 
 & \leq & 
2\, C_{m,p}\,|t|^{-p} \min\{1,|t|^2\}\, \frac{|z|^{m-1}}{\rho^{m-1}} \\ 
 & \leq &  
C_{m,p}\,|t|^{-p} \min\{1,|t|^2\}\, \left(\max\{|t|,|t|^3\}\right)^{m-1} |z|^{m-1}. 
\ene 
To simplify the dependence in $t$, note that in case $|t| \leq 1$, 
$$ 
|t|^{-p} \min\{1,|t|^2\}\, \left(\max\{|t|,|t|^3\}\right)^{m-1} = |t|^{(m+1)-p}, 
$$ 
while the left expression is equal to $|t|^{3(m-1)-p}$ in case $|t| \geq 1$. 
 
Also note that the condition $|tz| \leq 1$ is fulfilled automatically, 
as long as $|t^3 z| \leq 1$ and $|z| \leq 1$. Therefore, recalling also that  
$P_k = a_k$ for $k \leq m-2$, we obtain: 
 
\vskip5mm 
{\bf Proposition 8.3.} {\it If\, $0 < |z| \leq 1$ and $|t^3 z| \leq 1$,  
then for all $p = 0,1\dots,m$, 
$$ 
\frac{d^p}{dt^p}\, e^{W_z(t)} = \frac{d^p}{dt^p}\,  
\bigg(1 + \sum_{k=1}^{m-2} P_k(it)\,z^k\bigg) +  
A \left(|t|^{m+1-p} + |t|^{3(m-1)-p}\right) |z|^{m-1}, 
$$ 
where $|A| \leq C_{m,p} = (4^m (1+C))^p\ e^{2^m C}$. 
}


\vskip10mm 
\section{{\bf Edgeworth-type expansion for the functions $v(tz)^{1/z^2}$}} 
\setcounter{equation}{0} 
 
\vskip2mm 
Assume that $v(t)$ is $s$-times differentiable, $s \geq 2$, and  
not vanishing for $|t| \leq c$ $(c > 0)$, and such that $v(0)=1$, $v'(0) = 0$,  
$v''(0)=-1$. For  
$$ 
v_z(t) =  v(tz)^{1/z^2} 
$$  
define the approximating functions 
$$ 
u_m(t) = u_m(t,z) = e^{-t^2/2}\, \bigg(1 + \sum_{k=1}^{m-2} P_k(it)\, z^k\bigg), 
\qquad m = [s], 
$$ 
where the polynomials $P_k$ are based on the cumulants 
$\gamma_3,\dots,\gamma_m$ of $v$. Put $m'(s) = s-p$, 
$$ 
m''(s) = 3(m-2) + \max\{s+p, (s-1)p\}. 
$$ 
In particular, $m'(s) = s$ and $m''(s) = s + 3(m-2)$ in case $p=0$. 
Note that $m''(s) \leq 2m^2$ in all admissible cases. 
 
In this section, the relation (1.8) is established in the following more  
general form. 
 
\vskip5mm 
{\bf Proposition 9.1.} {\it Let $s \geq 3$. Given $z$ real, $0< |z| \leq 1$,  
in the interval $|t^3 z| \leq c^3$, for all $p = 0,1,\dots,m$, 
\be 
\big|v_z^{(p)}(t) - u_m^{(p)}(t)\big| \leq  
\big(|t|^{m'} + |t|^{m''}\big)\, e^{-t^2/2}\ |z|^{s-2}\, \ep(z), 
\en 
 
\vskip2mm 
\noindent 
where $\ep(z) \rightarrow 0$, as $z \rightarrow 0$. 
Moreover, if $s \geq 2$ and $v(t)$ has $(m+1)$ continuous derivatives, then 
with some constant $A$, and with $m'$, $m''$, corresponding to $s = m+1$, 
\be 
\big|v_z^{(p)}(t) - u_m^{(p)}(t)\big|  \leq  
A\, \big(|t|^{m'} + |t|^{m''}\big)\, e^{-t^2/2}\ |z|^{m-1}. 
\en 
 
} 
 
\vskip4mm 
We will refer to (9.1) and (9.2) as the scenarios 1 and 2, respectively.  
Note that in the second case, although $v$ has cumulants up to order $m+1$,  
we require that $\gamma_{m+1}$ does not  
participate in the definition of the polynomials $P_k$. In particular,  
the value $m=2$ is covered in (9.2), and then $u_m(t) = e^{-t^2/2}$  
(that is, $P_1$ is not present). 
 
\vskip2mm 
{\bf Proof.} Without loss of generality, assume $c=1$. Write 
$v_z(t) = e^{-t^2/2}\,w_z(t)\, e^{h_z(t)}$, where 
$$ 
w_z(t) = e^{W_z(t)}, \qquad 
W_z(t) = \sum_{k=3}^{m} \frac{\gamma_k}{k!}\, (it)^k z^{k-2}, 
$$ 
$$ 
\psi_z(t) = \frac{1}{2}\,t^2 + \frac{1}{z^2}\,\log v(tz) = 
\log\big(e^{t^2/2} v_z(t)\big), \qquad h_z(t) = \psi_z(t) - W_z(t). 
$$ 
By the definition of $a_k$ and $P_k$, 
$$ 
w_z(t) = 1 + \sum_{k=1}^{m-2} P_k(it)\,z^k + R_z(t), \qquad 
R_z(t) = \sum_{k=m-1}^\infty a_k(it)\,z^k. 
$$ 
Therefore, 
$$ 
v_z(t) =  u_m(t)\, e^{h_z(t)} + R_z(t)\, g(t)\, e^{h_z(t)}, \qquad 
g(t) = e^{-t^2/2}. 
$$ 
 
Given $p = 0,1\dots,m$, we differentiate this representation 
according to the Leibnitz rule: 
\begin{eqnarray} 
v_z^{(p)}(t) - u_m^{(p)}(t) \, = \,  I_1 + I_2 + I_3 
 & = & u_m^{(p)}(t)\, (e^{h_z(t)} - 1) \nonumber \\ 
 & & + \ 
\sum_{k=1}^{p} C_p^k\, u_m^{(p-k)}(t)\, (e^{h_z(t)})^{(k)} \nonumber \\ 
 & & \, + \, 
\sum_{k=0}^p C_p^k\, \left(R_z(t)\, g(t)\right)^{(k)}\, (e^{h_z(t)})^{(p-k)}, 
\end{eqnarray} 
where $C_p^k = \frac{p!}{k! (p-k)!}$ are the combinatorial coefficients. 
Note that when $p=0$, the second term $I_2$ is vanishing. 
 
\vskip2mm 
{\bf Estimation of} $I_1$. 
 
In Corollary 4.2 it is shown that, if $|z| \leq 1$ and $|t^3 z| \leq 1$,  
the functions $h_z$ and their derivatives are uniformly bounded and admit the bounds 
\be 
|h_z^{(p)}(t)| \leq |t|^{s-p} |z|^{s-2} \ep_p(z), \qquad p = 0,1,\dots,m, 
\en 
where each $\ep_p(z)$ is defined in $|z| \leq 1$  
and satisfies $\ep_p(z) \rightarrow 0$, as $z \rightarrow 0$.  
Moreover, if $v$ has $m+1$ continuous derivatives, then we have a sharper estimate 
\be 
|h_z^{(p)}(t)| \leq A_p\, |t|^{(m+1)-p}\, |z|^{m-1} 
\en 
with some constants $A_p$. In particular, when $p=0$, these bounds correspondingly  
give 
\be 
|e^{h_z(t)} - 1| \leq |t|^s |z|^{s-2} \ep(z), \qquad 
|e^{h_z(t)} - 1| \leq A_0\, |t|^{m+1}\, |z|^{m-1}, 
\en 
with some $\ep(z) \rightarrow 0$, as $z \rightarrow 0$, and a constant $A_0$.  
 
On the other hand, since every $P_k$ has degree $3k \leq 3(m-2)$, and $|z| \leq 1$, 
for all $p = 0,1,\dots,m$, $m \geq 3$, 
$$ 
\bigg|\frac{d^p}{dt^p}\,\bigg(1 + \sum_{k=1}^{m-2} P_k(it)\,z^k\bigg)\bigg| \leq  
C\big(1 + |t|^{3(m-2)-p}\big) 
$$ 
with some constant $C$, depending on $m$, $p$ and the cumulants $\gamma_k$'s.  
Since also 
\be 
|g^{(p)}(t)| = \bigg|\frac{d^p}{dt^p}\,e^{-t^2/2}\bigg| \leq C_p\, (1 + |t|^p)\,e^{-t^2/2}, 
\en 
we get, by the Leibnitz rule, 
\be 
|u_m^{(p)}(t)| \leq C\, \big(1 + |t|^{3(m-2) + p}\big)\,e^{-t^2/2}, 
\en 
where we allow the constants depend on $m$, $p$ and the cumulants  
$\gamma_3, \dots, \gamma_m$. For $m=2$, $u_m(t) = e^{-t^2/2}$, so (9.8) 
holds in this case as well ($p = 0,1,2$). 
Combining this with (9.6), we correspondingly arrive at 
\begin{eqnarray} 
|I_1| & \leq & \big(|t|^s + |t|^{s+p+3(m-2)}\big)\,e^{-t^2/2}\ |z|^{s-2} \ep(z),  
\\ 
|I_1| & \leq & A\, \big(|t|^{m+1} + |t|^{p+(4m-5)}\big)\,e^{-t^2/2}\ |z|^{m-1} 
\end{eqnarray} 
with some constant $A$ and $\ep(z) \rightarrow 0$, as $z \rightarrow 0$. 
As a result, we obtain the desired bounds on the first term $I_1$ in (9.3) 
for both scenarios.

\vskip2mm 
{\bf Estimation of} $I_2$. 
 
To treat the second term, more precisely -- the products 
$u_m^{(p-k)}(t)\, (e^{h_z(t)})^{(k)}$, assume that $p \geq 1$. 
By formula (2.3), for any $k = 1,\dots,p$, 
$$ 
(e^{h_z(t)})^{(k)} = e^{h_z(t)}\, k! \sum \prod_{r=1}^k \frac{1}{p_r!}  
\bigg(\frac{h_z^{(r)}(t)}{r!}\bigg)^{p_r}, 
$$ 
where the summation is performed over all non-negative integer 
solutions $(p_1,\dots,p_k)$ to the equation 
$p_1 + 2 p_2 + \dots + k p_k = k$. From (9.4)-(9.5) we get for the two scenarios 
$$ 
|h_z^{(r)}(t)|^{p_r} \leq |t|^{(s-r) p_r} |z|^{(s-2) p_r} \ep_r(z)^{p_r}, \qquad 
|h_z^{(r)}(t)|^{p_r} \leq A_r^{p_r}\, |t|^{((m+1)-r) p_r}\, |z|^{(m-1) p_r}. 
$$ 
After multiplication of these inequalities over all $r = 1,\dots,k$  
(separately in both scenarios), using $1 \leq p_1 + \dots + p_k \leq k$  
together with 
$$ 
s - k \leq \sum_{r=1}^k\, (s-r) p_r \leq (s-1)k, \qquad 
(m+1) - k \leq \sum_{r=1}^k\, ((m+1)-r) p_r \leq mk, 
$$ 
we obtain that 
\be 
\big|(e^{h_z(t)})^{(k)}\big| \leq \big(|t|^{s-k} + |t|^{(s-1)k}\big)  
|z|^{s-2} \ep(z), 
\en 
\be 
\big|(e^{h_z(t)})^{(k)}\big| \leq A\, \big(|t|^{(m+1)-k} + |t|^{mk}\big) |z|^{m-1} 
\en 
with some constant $A$ and $\ep(z) \rightarrow 0$, as $z \rightarrow 0$.  
One may combine these bounds with (9.8), which immediately yields 
$$ 
\big|u_m^{(p-k)}(t)\, (e^{h_z(t)})^{(k)}\big| \leq  
\big(|t|^{s-k} + |t|^{(s-1)k + (p-k) + 3(m-2)}\big)\,e^{-t^2/2}\ |z|^{s-2}\, \ep(z), 
$$ 
$$ 
\big|u_m^{(p-k)}(t)\, (e^{h_z(t)})^{(k)}\big| \leq A\, 
\big(|t|^{(m+1)-k} + |t|^{mk + (p-k) + 3(m-2)}\big)\,e^{-t^2/2}\ |z|^{m-1}. 
$$ 
Since $k$ varies from 1 to $p$, the right-hand sides can be made independent of $k$, 
and we arrive at the desired bounds on the second term $I_2$ in (9.3), 
needed for the values $p \geq 1$: 
\be 
|I_2| \leq  
\big(|t|^{s-p} + |t|^{(s-1) p + 3(m-2)}\big)\,e^{-t^2/2}\ |z|^{s-2} \ep(z), 
\en 
\be 
|I_2| \leq  
\big(|t|^{(m+1)-p} + |t|^{m p + 3(m-2)}\big)\,e^{-t^2/2}\ |z|^{m-1}. 
\en

\vskip2mm 
{\bf Estimation of} $I_3$. 
 
Now, let us turn to the third term, i.e., to the products 
$\left(R_z(t)\, g(t)\right)^{(k)}\, (e^{h_z(t)})^{(p-k)}$. 
By Proposition 8.3, for all $p = 0,1,\dots,m$, 
$$ 
\big|R_z^{(p)}(t)\big| \leq C \big(|t|^{(m+1)-p} + |t|^{3(m-1)-p}\big) |z|^{m-1}. 
$$ 
Using (9.7) and the Leibnitz formula, the latter gives, for all $k = 0,\dots,p$, 
\be 
\big|\big(R_z(t)\, g(t)\big)^{(k)}\big| \leq  
C \big(|t|^{(m+1)-k} + |t|^{3(m-1) + k}\big)\,e^{-t^2/2}\ |z|^{m-1}. 
\en 
 
{\bf Case} $p=0$. Then necessarily $k=0$, and the above inequality yields 
\be 
|I_3| \leq C \big(|t|^{m+1} + |t|^{3(m-1)}\big)\,e^{-t^2/2}\ |z|^{m-1}. 
\en 
It has only to be compared with (9.9)-(9.10). In the second scenario, one clearly 
obtains from (9.10) and (9.16) that 
$$ 
|v_z(t) - u_m(t)| \leq |I_1| + |I_3| \leq 
C \big(|t|^{m+1} + |t|^{5m-5}\big)\,e^{-t^2/2}\ |z|^{m-1}. 
$$ 
This proves (9.2) in case $p=0$.  
 
In the first scenario, in (9.16) just write 
$|z|^{m-1} = |z|^{s-2}\, \widetilde\ep(z)$ with $\widetilde\ep(z) \rightarrow 0$,  
as $z \rightarrow 0$. Together with (9.9) this leads to a similar estimate 
$$ 
|v_z(t) - u_m(t)| \leq |I_1| + |I_3| \leq 
C \big(|t|^s + |t|^{s + 3(m-2)}\big)\,e^{-t^2/2}\ |z|^{s-2}, 
$$ 
proving (9.1) in case $p=0$. Thus, Proposition 9.1 is proved in this case.

\vskip5mm 
{\bf Case} $1 \leq p \leq m$. If $k = p$, the absolute value of the product 
$$ 
\left(R_z(t)\, g(t)\right)^{(k)}\, (e^{h_z(t)})^{(p-k)} = 
\left(R_z(t)\, g(t)\right)^{(p)}\, e^{h_z(t)} 
$$ 
may be estimated according to (9.15) by 
\be 
C \big(|t|^{(m+1)-p} + |t|^{3(m-1) + p}\big)\,e^{-t^2/2}\ |z|^{m-1}. 
\en 
As in the previous step, $|z|^{m-1}$ may be replaced here with $|z|^{s-2}\, \ep(z)$. 
 
If $0 \leq k \leq p-1$, by (9.11)-(9.12) for the two scenarios we have 
\be 
\big|(e^{h_z(t)})^{(p-k)}\big| \leq \big(|t|^{s-(p-k)} + |t|^{(s-1)(p-k)}\big)  
|z|^{s-2} \ep(z), 
\en 
\be 
\big|(e^{h_z(t)})^{(p-k)}\big| \leq A\, \big(|t|^{(m+1)-(p-k)} + |t|^{m (p-k)}\big) |z|^{m-1}. 
\en 
It remains to multiply these inequalities by (9.15). 
In this step we consider the two scenarios separately. 
 
\vskip5mm 
{\bf Scenario} 1 (The inequality (9.1)): When multiplying (9.15) by (9.18)  
and looking for the maximal power of $|t|$, notice that 
$$ 
(3(m-1) + k) + (s-1)(p-k) 
$$ 
is maximized for $k=0$, and for this value it is equal to $3(m-1) + (s-1)p$. 
Therefore, 
$$ 
\big|\big(R_z(t)\, g(t)\big)^{(k)}\, (e^{h_z(t)})^{(p-k)}\big| \leq 
\big(|t|^{(m+1) + (s-p)} + |t|^{3(m-1) + (s-1)p}\big)\,e^{-t^2/2}\  
|z|^{(m-1) + (s-2)} \ep(z). 
$$ 
 
Now, comparing with (9.9), (9.13) and (9.17), we see that the smallest power of  
$|t|$ in these inequalities is $m' = s-p$. Hence, we do not loose much by writing 
$$ 
\big|\big(R_z(t)\, g(t)\big)^{(k)}\, (e^{h_z(t)})^{(p-k)}\big| \leq 
\big(|t|^{m'} + |t|^{3(m-1) + (s-1)p}\big)\,e^{-t^2/2}\ |z|^{(m-1) + (s-2)} \ep(z), 
$$ 
which holds for all $k = 0,1,\dots,p$. 
To simplify, let us note that $|t|^{3(m-1)} |z|^{(m-1)} \leq 1$, which leads to 
\be 
\big|\big(R_z(t)\, g(t)\big)^{(k)}\, (e^{h_z(t)})^{(p-k)}\big| \leq 
\big(|t|^{m'} + |t|^{(s-1)p}\big)\,e^{-t^2/2}\ |z|^{s-2} \ep(z). 
\en 
In addition, the largest power of $|t|$ in (9.9), (9.13), (9.17) and (9.20) is  
$$ 
m'' = 3(m-2) + \max\{s+p, (s-1)p\}. 
$$ 
Hence, $I_3 \leq \big(|t|^{m'} + |t|^{m''}\big)\,e^{-t^2/2}\ |z|^{s-2} \ep(z)$ 
and similarly for $I_1$ and $I_2$. This proves (9.1).

\vskip5mm 
{\bf Scenario} 2 (The inequality (9.2)): This case can be dealt with along 
the lines of scenario~1 by letting $s \rightarrow m+1$. Or, repeating  
the previous arguments, note that when multiplying (9.15) by (9.19),  
the expression $3(m-1) + k + m(p-k)$  
is maximized for $k=0$, and for this value it is equal to $3(m-1) + mp$. Therefore, 
$$ 
\big|\big(R_z(t)\, g(t)\big)^{(k)}\, (e^{h_z(t)})^{(p-k)}\big| \leq A\, 
\big(|t|^{2(m+1) - p} + |t|^{3(m-1) + mp}\big)\,e^{-t^2/2}\  
|z|^{2(m-1)}. 
$$ 
In (9.10), (9.14) and (9.17) the smallest power of $|t|$ is  
$m' = (m+1)-p$. Hence, 
$$ 
\big|\big(R_z(t)\, g(t)\big)^{(k)}\, (e^{h_z(t)})^{(p-k)}\big| \leq A\, 
\big(|t|^{m'} + |t|^{3(m-1) + mp}\big)\,e^{-t^2/2}\ |z|^{2(m-1)}. 
$$ 
Again using $|t|^{3(m-1)} |z|^{m-1} \leq 1$, we get  
\be 
\big|\big(R_z(t)\, g(t)\big)^{(k)}\, (e^{h_z(t)})^{(p-k)}\big| \leq 
\big(|t|^{m'} + |t|^{mp}\big)\,e^{-t^2/2}\ |z|^{m-1}. 
\en 
In addition, the largest power of $|t|$ in (9.10), (9.14), (9.17) and (9.21) is 
$$ 
m'' = 3(m-2) + \max\{m+p+1, mp\}. 
$$ 
Hence, $I_3 \leq A \big(|t|^{m'} + |t|^{m''}\big)\,e^{-t^2/2}\ |z|^{m-1}$ 
and similarly for $I_1$ and $I_2$ in case $p \geq 1$. 
 
This proves (9.2) and Proposition 9.1.

\vskip5mm 
One may unite Proposition 9.1 (first part) with Propositions 5.1-5.2  
for the case $2 \leq s < 3$, if we do not care about  
polynomial factors in front of $e^{-t^2/2}$.

\vskip5mm 
{\bf Corollary 9.2.} {\it There is a function $T_z \rightarrow +\infty$,  
as $z \rightarrow 0$ $(0 \leq |z| \leq 1)$, such that in the interval  
$|t| \leq T_z$, for all $p = 0,1,\dots,m$, $m = [s]$, 
$$ 
\big|v_z^{(p)}(t) - u_m^{(p)}(t) \big| \leq \ep(z) |z|^{s-2}\,e^{-t^2/4} 
$$ 
with $\ep(z) \rightarrow 0$. Moreover, up to some constant $c>0$, one can choose  
$T_z = c\, |z|^{-1/3}$ in case $s \geq 3$ and $T_z = c\, |z|^{-(s-2)/s}$ in case  
$2 < s < 3$. 
}


\vskip10mm 
\section{{\bf Proof of Theorem 1.3}} 
\setcounter{equation}{0} 
 
\vskip2mm 
Again, let $v(t)$ be $s$-times differentiable, $s \geq 2$, such that  
$v(0)=1$, $v'(0) = 0$, $v''(0)=-1$. Note that $v$ is not vanishing in some 
interval, containing the origin.  
 
Let us return to the family of the functions 
$$ 
u_m(t,z) =  
e^{-t^2/2}\, \bigg(1 + \sum_{k=1}^{m-2} P_k(it) z^k\bigg), 
$$ 
where $m = [s]$, and the polynomials $P_k$ are based on the cumulants 
$\gamma_3,\dots,\gamma_m$ of $v$. 
In order to approximate the powers $v_n(t) =  v(\frac{t}{\sqrt{n}})^n$, 
one uses the values $z = 1/\sqrt{n}$, leading to the the approximating functions 
$$ 
u_m(t) \, = \, u_m(t,n^{-1/2}) \, = \, 
e^{-t^2/2}\, \bigg(1 + \sum_{k=1}^{m-2} P_k(it)\, n^{-k/2}\bigg). 
$$ 
On the other hand, when $z=1$, we deal with the 
projection operators $T_m$, i.e., with the functions 
$$ 
e_m(t) = u_m(t,1) = e^{-t^2/2}\, \bigg(1 + \sum_{k=1}^{m-2} P_k(it)\bigg). 
$$ 
Theorem 1.3 is a particular case of the following more general proposition. 
 
\vskip5mm 
{\bf Proposition 10.1.} {\it For all $p = 0,1\dots,m$, 
and all $|t| \leq cn^{1/6}$, 
\be 
\frac{d^p}{dt^p}\, (v_n(t) - u_m(t)) = n\, \frac{d^p}{dt^p} 
\left[\bigg(v\big(\frac{t}{\sqrt{n}}\big) - e_m\big(\frac{t}{\sqrt{n}}\big)\bigg)\,  
e^{-t^2/2}\right] + r_n 
\en 
with 
\be 
|r_n| \leq \big(1 + |t|^{4m^2}\big)\, e^{-t^2/2}\  
\bigg(\frac{C}{n^{(m-1)/2}} + \frac{\ep_n}{n^{s-2}}\bigg). 
\en 
Here $C$, $c$ and $\ep_n$ are some positive constants, such that $\ep_n \rightarrow 0$, 
as $n \rightarrow \infty$. 
} 
 
\vskip5mm 
It is worthwile noting that Proposition 9.1 and thus relation (1.8) can be  
obtained on the basis of (10.1)-(10.2) as well, using the property that  
$v(t)$ and $e_m(t)$ have equal derivatives up to order $m$ and both are  
$s$-times differentiable.  
However, we have chosen a different road of proof and will derive Proposition 10.1  
by virtue of Proposition 9.1 (its second part). 
 
In order to show how it applies, apply the binomial formula to obtain that 
\bee 
v_n(t) - u_m(t)\, =\, \sigma_{n1} + \sigma_{n2} + \sigma_{n3} 
 & = & 
e_m\big(\frac{t}{\sqrt{n}}\big)^n - u_m(t) \\ 
 & & +\ 
n \left[v\big(\frac{t}{\sqrt{n}}\big) - e_m\big(\frac{t}{\sqrt{n}}\big)\right] 
e_m\big(\frac{t}{\sqrt{n}}\big)^{n-1} 
 \\ 
 & & +\, 
\sum_{k=2}^n C_n^k  
\left[(v\big(\frac{t}{\sqrt{n}}\big) - e_m\big(\frac{t}{\sqrt{n}}\big)\right]^k 
e_m\big(\frac{t}{\sqrt{n}}\big)^{n-k}. 
\ene 
 
Thus, $\sigma_{n2}$ is almost the term which appears on the right-hand  
side of (10.1), provided that $e_m(\frac{t}{\sqrt{n}})^{n-1}$  
is replaced with the characteristic function of the standard normal law. 
 
The first term $\sigma_{n1} = e_m(\frac{t}{\sqrt{n}})^n - u_m(t)$ is of  
the same nature as $v_n(t) - u_m(t)$, assuming that $e_m$ plays the role  
of $v$. At this point, let us recall that, by Proposition 7.4,  
$T_m e_m = e_m$, and moreover, that $e_m$ generates the same polynomials  
$P_k$ as $v$. Hence, Proposition 9.1, being applied to $e_m$  
in place of $v$ with $z = 1/\sqrt{n}$, provides the bound on  
the derivatives of $e_m(\frac{t}{\sqrt{n}})^n - u_m(t)$. 
Since $e_m$ is analytic, the second assertion (9.2) of Proposition 9.1  
is more accurate. Namely, if $e_m(t)$ is not vanishing in the interval  
$|t| \leq c$ (which is true 
with some constant $c>0$, depending on the cumulants, only), it gives: 
 
\vskip5mm 
{\bf Lemma 10.2.} {\it For all $p = 0,1,\dots,m$ and all $|t| \leq c n^{1/6}$,  
\be 
\big|\sigma_{n1}^{(p)}(t)\big|  \leq  
A \big(1 + |t|^{2m^2}\big)\, e^{-t^2/2}\ n^{-(m-1)/2}, 
\en 
where $c$ and $A$ are some positive constants, depending on the cumulants 
$\gamma_3,\dots,\gamma_m$. 
} 
 
\vskip5mm 
Moreover, using a similar argument one may estimate the derivatives of the  
functions $e_m(t)^k$, which appear both in $\sigma_{n2}$ and $\sigma_{n3}$. 
Apply (9.2) with $v = e_m$ and $z = 1/\sqrt{k}$ to get 
\be 
\frac{d^p}{dt^p}\ e_m\big(\frac{t}{\sqrt{k}}\big)^k = \frac{d^p}{dt^p}\ 
e^{-t^2/2}\, \bigg(1 + \sum_{j=1}^{m-2} \frac{P_j(it)}{k^{j/2}}\bigg) + 
A \big(1 + |t|^{2m^2}\big) e^{-t^2/2}\ k^{-(m-1)/2}, 
\en 
where $A = A_k(t)$ is a bounded quantity in the interval $|t| \leq ck^{1/6}$. 
Putting $\alpha = \sqrt{\frac{k}{n}}$ and replacing the variable $t$ with $\alpha t$,  
we obtain 
\be 
\frac{d^p}{dt^p}\ e_m\big(\frac{t}{\sqrt{n}}\big)^k = \frac{d^p}{dt^p}\ 
e^{-\alpha^2 t^2/2}\, \bigg(1 + \sum_{j=1}^{m-2} \frac{P_j(i\alpha t)}{k^{j/2}}\bigg) + 
B \big(1 + |t|^{2m^2}\big) e^{-\alpha^2 t^2/2}, 
\en 
where now $B = B_k(t)$ is bounded in $|t| \leq cn^{1/6}$. Every $P_j$ is  
a polynomial of degree at most $3j \leq 3m$, so all its derivatives of  
order up to $m$ can be bounded by $C(1 + |t|)^{3m}$ on the whole real line.  
Hence, using $\alpha \leq 1$ and $k \geq 1$, the same is true for  
the polynomial in the large bracket of (10.5). Using the Leibnitz rule,  
it then follows from (10.5) that:

\vskip5mm 
{\bf Lemma 10.3.} {\it For all $p = 0,1,\dots,m$, $k = 0,1,\dots,n$, in the interval  
$|t| \leq c n^{1/6}$ 
$$ 
\bigg|\frac{d^p}{dt^p}\ e_m\big(\frac{t}{\sqrt{n}}\big)^k\bigg| \leq 
C \big(1 + |t|^{2m^2}\big) e^{-kt^2/(2n)} 
$$ 
with some positive constants $c$ and $C$. 
} 
 
\vskip5mm 
The particular case $k = n-1$ in (10.4) should be investigated in more  
detail by replacing (10.5) with a more accurate relation, namely 
$$ 
\frac{d^p}{dt^p}\ e_m\big(\frac{t}{\sqrt{n}}\big)^{n-1} = \frac{d^p}{dt^p}\ 
e^{-\alpha^2 t^2/2}\, \bigg(1 + \frac{1}{\sqrt{n-1}}\,\sum_{j=1}^{m-2}  
\frac{P_j(i\alpha t)}{(n-1)^{(j-1)/2}}\bigg) +  
\frac{B\big(1 + |t|^{2m^2}\big)}{\sqrt{n-1}}\, e^{-\alpha^2 t^2/2} 
$$ 
(assuming $n \geq 2$). Repeating the same argument concerning the growth  
of the polynomials $P_j$ and its derivatives, and noting that, for  
$\alpha = \sqrt{\frac{n-1}{n}}$, we have 
$$ 
\bigg|\frac{d^p}{dt^p}\, e^{-\alpha^2 t^2/2} - \frac{d^p}{dt^p}\, e^{-t^2/2}\bigg| 
\leq \frac{C}{n}\,\big(1 + |t|^{p+2}\big)\,e^{-t^2/2}, \qquad  
|t| \leq \sqrt{n}, 
$$ 
we arrive at the following bound (which also holds in the missing case $n=1$): 
 
\vskip5mm 
{\bf Lemma 10.4.} {\it For all $p = 0,1,\dots,m$, in the interval $|t| \leq c n^{1/6}$ 
$$ 
\bigg|\frac{d^p}{dt^p}\ e_m\big(\frac{t}{\sqrt{n}}\big)^{n-1} -  
\frac{d^p}{dt^p}\, e^{-t^2/2}\bigg|  
\leq \frac{C}{\sqrt{n}}\, \big(1 + |t|^{2m^2}\big)\, e^{-t^2/2} 
$$ 
with some positive constants $c$ and $C$. 
}

\vskip5mm 
Finally, let us bound the derivatives of $y(t) = v(t)-e_m(t)$ and of its  
powers, which appear both in $\sigma_{n2}$ and $\sigma_{n3}$ as well.  
To this aim, one may appeal to Corollary 7.5, giving, as $t \rightarrow 0$, 
\be 
y^{(r)}(t) = o\big(|t|^{s-r}\big), \quad {\rm for \ any} \ \ r = 0,\dots,m.  
\en 
In particular, 
$y(t)^k = o(|t|^{sk})$, for any $k \geq 1$. If $p \geq 1$, by the chain rule  
(cf. (2.3)), the $p$-th derivative of $y(t)^k$ represents a linear combination  
of the terms 
$$ 
b(t) = y(t)^{k - (k_1 + \dots + k_p)}\, (y'(t))^{k_1} \dots (y^{(p)}(t))^{k_p} 
$$ 
over all integer tuples $(k_1,\dots,k_p)$, such that  
$k_1 + 2k_2 + \dots + p k_p = p$ and  
$k_1 + \dots + k_p \leq k$ ($k_j \geq 0$). By (10.6), we have  
$b(t) = o(|t|^S)$,  
where 
$$ 
S = s\big(k - (k_1 + \dots + k_p)\big) + \sum_{r=1}^p (s-r) k_r = sk-p. 
$$ 
Hence, $\frac{d^p}{dt^p}\, \big(v(t) - e_m(t)\big)^k = o\big(|t|^{sk-p}\big)$. 
Since $sk-p \leq sm$ for $1 \leq k \leq m$, we obtain:

\vskip5mm 
{\bf Lemma 10.5.} {\it Let $0 < \alpha < \frac{1}{2}$ and $c>0$ be given. 
For some $\ep_n \rightarrow 0$, for all $p = 0,1,\dots,m$ and $k = 1,\dots,m$,  
we have, uniformly in the interval $|t| \leq c n^{\alpha}$, 
$$ 
\bigg|\frac{d^p}{dt^p}\, \bigg(v\big(\frac{t}{\sqrt{n}}\big) -  
e_m\big(\frac{t}{\sqrt{n}}\big)\bigg)^k\bigg| \leq 
\ep_n \big(1 + |t|^{sm}\big)\, n^{-sk/2}. 
$$ 
} 
 
\vskip2mm 
{\bf Proof of Proposition 10.1.} Using Lemmas 10.4 and 10.5 (with $k=1$),  
we see that in $\sigma_{n2}$ one may replace the term 
$e_m(\frac{t}{\sqrt{n}})^{n-1}$  
with $e^{-t^2/2}$ at the expense of an error, not exceeding 
\be 
n \cdot \ep_n \big(1 + |t|^{sm}\big)\, n^{-s/2} \cdot 
\frac{C}{\sqrt{n}}\, \big(1 + |t|^{2m^2}\big)\, e^{-t^2/2} \leq 
\frac{\ep_n'}{n^{(s-1)/2}}\, \big(1 + |t|^{4m^2}\big)\, e^{-t^2/2}, 
\en 
where $\ep_n' \rightarrow 0$. The same is true for the first $m$  
derivatives of $\sigma_{n2}$. 
 
Now consider the products  
$y_k(t) = \big((v(\frac{t}{\sqrt{n}}) - e_m(\frac{t}{\sqrt{n}})\big)^k\,  
e_m\big(\frac{t}{\sqrt{n}}\big)^{n-k}$ 
appearing in $\sigma_{n3}$ with $2 \leq k \leq n$. Writing 
$$ 
y_k^{(p)}(t) = \sum_{j=0}^{p} C_p^j \ \frac{d^j}{dt^j}\,  
\bigg(v\big(\frac{t}{\sqrt{n}}\big) - e_m\big(\frac{t}{\sqrt{n}}\big)\bigg)^k\, 
\frac{d^{p-j}}{dt^{p-j}}\,e_m\big(\frac{t}{\sqrt{n}}\big)^{n-k} 
$$ 
and combining Lemmas 10.3 and 10.5 (which give estimates that are  
independent of $j$), we get 
\bee 
\big|y_k^{(p)}(t)\big| & \leq &  
2^p \cdot \ep_n \big(1 + |t|^{sm}\big)\, n^{-sk/2} \cdot 
C \big(1 + |t|^{2m^2}\big)\, e^{-(n-k) t^2/(2n)} \\ 
 & \leq &  
\frac{\ep_n''}{n^{sk/2}}\, \big(1 + |t|^{4m^2}\big)\, e^{-(n-k) t^2/(2n)}. 
\ene 
Therefore, 
\begin{eqnarray} 
|\sigma_{n3}^{(p)}(t)| & \leq &  
\sum_{k=2}^{n} C_n^k \, \big|y_k^{(p)}(t)\big| \nonumber \\ 
 & \leq &  
\ep_n''\, \big(1 + |t|^{4m^2}\big)\sum_{k=2}^{n}  
C_n^k \frac{1}{n^{sk/2}}\, \,e^{-(n-k) t^2/(2n)} \nonumber \\ 
 & = & 
\ep_n''\, \big(1 + |t|^{4m^2}\big)\, e^{-t^2/2}\, 
\bigg(\big(1 + n^{-s/2}\, e^{t^2/2n}\big)^n - 1 - n^{-(s-2)/2}\, e^{t^2/2n}\bigg). 
\end{eqnarray} 
For $s > 2$, we have $\delta_n = n^{-s/2}\, e^{t^2/2n} = o(1/n)$  
uniformly in the interval $|t| \leq cn^{1/6}$. So, 
$$ 
(1 + \delta_n)^n =  
e^{n \log(1 + \delta_n)} = e^{n (\delta_n + O(\delta_n^2))} = 
1 + n\delta_n + \frac{1}{2}\,(n\delta_n)^2 + n O\big(\delta_n^2\big). 
$$ 
Hence, for all $n$ large enough, the expression in the large bracket in  
(10.8) does not exceed 
$$ 
\frac{1}{2}\,(n\delta_n)^2 + O\big(\delta_n^2\big) \leq \frac{1}{n^{s-2}} +  
O\bigg(\frac{1}{n^{s-1}}\bigg). 
$$ 
It remains to compare this bound with (10.7) and (10.3), and then we arrive  
at (10.2). 
 
Finally, in the case $s=2$, the expression in the large brackets in (10.8)  
is uniformly bounded in $|t| \leq cn^{1/6}$. 
Thus Proposition 10.1 is proved.


\vskip10mm 
\section{{\bf Liouville Fractional Integrals and Derivatives}} 
\setcounter{equation}{0} 
 
\vskip2mm 
In this section we recall basic definitions and some results on Liouville fractional integrals and derivatives, and refer to [S-K-M], [K-S-T]  
for proofs and a more detailed exposition.  
At the end of the section we also formulate some special estimates 
for such operators. The proof of Proposition 11.3 below is rather routine  
and is therefore postponed to the next section. 
 
Let $\alpha$ denote a real number with $0 < \alpha < 1$, 
and let $y = y(t)$ denote a (measurable) function defined for $t > 0$. 
The Liouville left- and right-sided fractional integrals on the positive half-axis 
$\R^+ = (0,+\infty)$ of order $\alpha$ are defined by 
$$ 
(I^\alpha_{0+} y)(x) = \frac{1}{\Gamma(\alpha)} \int_0^x  
\frac{y(t)\,dt}{(x-t)^{1-\alpha}}, 
\qquad 
(I^\alpha_{-} y)(x) = \frac{1}{\Gamma(\alpha)} \int_x^{+\infty}  
\frac{y(t)\,dt}{(t-x)^{1-\alpha}} \qquad (x>0). 
$$ 
 
The equalities are understood in the usual way (as Lebesgues integrals),  
if $y$ is sufficiently "nice". According to a theorem by 
Hardy abd Littlewood, $I^\alpha_{0+}$ and $I^\alpha_{-}$ are extended  
and act as bounded linear operators from $L^p(\R^+)$  
to $L^q(\R^+)$, where $1 \leq p,q \leq +\infty$, if and only if 
$p < \frac{1}{\alpha}$ and $q = \frac{p}{1 - \alpha p}$. 
They represent particular cases of the so-called Liouville (or Riemann-Liouville) 
fractional calculus operators. 
 
The Liouville left- and right-sided fractional derivatives on the positive  
half-axis are defined by 
$$ 
(D^\alpha_{0+} y)(x) = \frac{d}{dx}\, (I^{1-\alpha}_{0+} y)(x), 
\qquad 
(D^\alpha_{-} y)(x) = \frac{d}{dx}\, (I^{1-\alpha}_{-} y)(x) 
\qquad (x>0). 
$$ 
The equalities are valid for sufficiently "nice" functions, including 
the class $C_0^\infty(\R^+)$ of all infintely differentiable functions 
on $\R^+$ with a compact support (which can be used to approximate functions 
from larger spaces). 
 
For example, for any complex number $\lambda$, such that  
${\rm Re}(\lambda) > 0$, 
\be 
(I^\alpha_{-} e^{-\lambda t})(x) = \lambda^{-\alpha}\, e^{-\lambda x}, \qquad 
(D^\alpha_{-} e^{-\lambda t})(x) = \lambda^{\alpha}\, e^{-\lambda x}, 
\en 
where the principal value of the power functions is used. 
 
We cite two standard facts about these operators (see [K-S-T], p.75 and p.83). 
 
\vskip5mm 
{\bf Proposition 11.1.} {\it For all sufficiently "good" functions $y$ on $\R^+$, 
$$ 
(D^\alpha_{0+} I^{\alpha}_{0+} y)(x) = y(x), \qquad 
(D^\alpha_{-} I^{\alpha}_{-} y)(x) = y(x). 
$$ 
The equalities are extended to the space $L^1(\R^+)$. Moreover, if  
additionally $y(x) = o(x^\alpha)$ for $x \rightarrow 0$, then 
$$ 
(I^{\alpha}_{0+} D^\alpha_{0+} y)(x) = y(x). 
$$ 
} 
 
Define the linear spaces $I^{\alpha}_{0+}(L^p(\R^+))$ and  
$I^{\alpha}_{-}(L^p(\R^+))$ as the images of $L^p(\R^+)$ under  
the operators $I^{\alpha}_{0+}$ and $I^{\alpha}_{-}$, respectively. 
 
\vskip5mm 
{\bf Proposition 11.2.} {\it For all sufficiently "good" functions $f$ and 
$g$ on $\R^+$, 
\be 
\int_0^{+\infty} f(x)\, (D_{0+}^\alpha g)(x)\,dx =  
\int_0^{+\infty} g(x)\, (D_{-}^\alpha f)(x)\,dx. 
\en 
The equality may be extended to all $f \in I^{\alpha}_{-}(L^p(\R^+))$  
and $g \in I^{\alpha}_{0+}(L^q(\R^+))$ with $p,q>1$, such that  
$\frac{1}{p} + \frac{1}{q} = 1+\alpha$. 
} 
 
\vskip5mm 
This is a formula for fractional integration by parts. 
 
Now, let $V$ be a function of bounded variation on the real line, also viewed 
as a finite measure, and denote by $|V|$ its variation (as a measure). 
Define the Fourier-Stieltjes transform 
$$ 
\hat V(x) = \int_{-\infty}^{+\infty} e^{itx}\,dV(t) \quad x \in \R. 
$$ 
 
For our purposes, the following proposition will play a crucial role in the study 
of the local limit theorem with fractional moments. 
 
\vskip5mm 
{\bf Proposition 11.3.} {\it Let $g(x) = \hat V(x)\, h(x)$, where  
$h(x)$ is a continuously differentiable function on the real line, such that,  
for a given integer $m \geq 0$ and $0 < \alpha < 1$, as $|x| \rightarrow \infty$, 
$$ 
|h(x)| + |h'(x)| = O\big(|x|^{-(2+m+\alpha)}\big). 
$$ 
If $\int_{-\infty}^{+\infty} |t|^{m+\alpha}\,d|V|(t) < +\infty$ and  
$V^{(k)}(0)=0$, for all $k = 0,\dots,m$, then $(D_{0+}^\alpha g)(x)$  
exists for all $x>0$ and satisfies with some constant $C$, independent of $x$, 
\be 
\big|(D_{0+}^\alpha g)(x)\big| \leq \frac{C}{(1+x)^\alpha}\, \int_{-\infty}^{+\infty} 
\min\big\{|u|,|u|^{\alpha'}\big\}\,|u|^m\,|V|(du). 
\en 
Here $\alpha' = \alpha$ in case $m = 0$, and $\alpha' = 0$ in case $m \geq 1$.  
In addition, for all $t$ real, 
\be 
\int_0^{+\infty} e^{it x}\, (D_{0+}^\alpha g)(x)\,dx =  
(-it)^\alpha 
\int_0^{+\infty}  e^{it x}\, g(x)\,dx. 
\en 
}

\vskip2mm 
More precisely, the Gaussian function $h(x) = e^{-x^2/2}$ and its derivatives 
will only be needed in this proposition to obtain the desired decay  
for the inverse Fourier-Stieltjes transform. 
 
\vskip5mm 
{\bf Proposition 11.4.} {\it For all functions $V$ and $h$ as in Proposition $11.3$,  
for all $t$ real, 
\be 
\bigg|\int_{-\infty}^{+\infty}  e^{it x}\, \hat V(zx)\,h(x)\,dx\bigg|\, \leq\, 
\frac{|z|^{m+\alpha}}{(1 + |t|)^\alpha}\, \ep(z), 
\en 
where $\ep(z)$ is bounded in $|z| \leq 1$ and satisfies  
$\ep(z) \rightarrow 0$, as $z \rightarrow 0$. 
}

\vskip5mm 
{\bf Proof.} Let $0 < z \leq 1$. We apply Proposition 11.3 with the function 
$V_z(u) = V(u/z)$ in place of $V$ in which case $\hat V_z(x) = \hat V(zx)$. 
Then, for the function $g_z(x) = \hat V(zx)\,h(x)$, the fractional 
derivative $(D_{0+}^\alpha g_z)(x)$ exists for all $x>0$ and satisfies (11.3). 
 
In order to unite both cases, use $|u|^{\alpha'} \leq |u|^{\alpha}$ 
for $|u| \geq 1$ and write (11.3) in a slightly weaker form 
\be 
\big|(D_{0+}^\alpha g_z)(x)\big| \leq \frac{Cz^m}{(1+x)^\alpha}\,\delta(z), 
\en 
where 
$$ 
\delta(z) = \int_{-\infty}^{+\infty} 
\min\big\{|zu|,|zu|^{\alpha}\big\}\,|u|^m\,|V|(du). 
$$ 
This integral is finite and behaves like $o(z^\alpha)$, as $z \rightarrow 0$. 
Indeed, split it into the two integrals in terms of the (finite positive)  
measure $W(du) = |u|^m\,|V|(du)$ as 
\be 
\delta(z) = z^\alpha I_0(z) + z^\alpha I_1(z) =  
z^\alpha \int_{|u| \leq 1/z} z^{1-\alpha} |u|\, dW(u) + 
z^\alpha \int_{|u| > 1/z} |u|^\alpha\, dW(u). 
\en 
By the moment assumption on $V$, we have $\int |u|^\alpha\, dW(u) < +\infty$, so 
$I_1(z) \rightarrow 0$, as $z \rightarrow 0$. As for the first integral, 
note that $z^{1-\alpha} |u| \leq |u|^\alpha$ in the region $|u| \leq 1/z$. 
Hence, the functions $f_z(u) = z^{1-\alpha} |u|\, 1_{\{|u| \leq 1/z\}}$ 
have an integrable majorant $f(u) = |u|^\alpha$ with respect to $W$.  
Since also $f_z(u) \rightarrow 0$, as $z \rightarrow 0$,  
one may apply the Lebesgue dominated convergence theorem, which gives 
$I_0(z) = \int f_z\, dW \rightarrow 0$. Thus, from (11.6)-(11.7), 
\be 
\big|(D_{0+}^\alpha g_z)(x)\big| \leq \frac{Cz^{m+\alpha}}{(1+x)^\alpha}\,\ep(z), 
\en 
where 
$\ep(z) = \frac{\delta(z)}{z^\alpha} \rightarrow 0$, as $z \rightarrow 0$, and  
$\sup_{0<z\leq 1} \ep(z) < +\infty$. 
 
\vskip2mm 
Now, using the bound (11.8) in (11.4), we get 
$|\int_0^{+\infty}  e^{it x}\, g_z(x)\,dx| \leq \frac{z^\alpha}{|t|^\alpha}\, \ep(z)$. 
Obviously, a similar inequality will hold as well when integrating 
over the negative half-axis. Therefore, 
$$ 
\bigg|\int_{-\infty}^{+\infty}  e^{it x}\, g_z(x)\,dx\bigg|\, \leq\, 
\frac{z^\alpha}{|t|^\alpha}\, \ep(z). 
$$ 
This estimate implies (11.5) in case of large values of $|t|$, say, when  
$|t| \geq 1$. The remaining range $|t| \leq 1$ can be treated by  
straightforward arguments. 
 
By the assumption on the decay of $h$, its Fourier transform 
$\hat h(t) = \int_{-\infty}^{+\infty}  e^{itx}\, h(x)\,dx$ is well-defined,  
bounded, and has bounded continuous derivatives up to order $m+1$.  
Introduce Taylor's approximation for $\hat h$ up to order $m$ at a given point  
$t$, i.e., the function 
$$ 
(S_m \hat h)(t,u) = \sum_{k=0}^m \frac{\hat h^{(k)}(t)}{k!}\, u^k, \qquad  
u \in \R. 
$$ 
From Taylor's theorem it follows that 
\be 
\big|\hat h(t + u) - (S_m \hat h)(t,u)\big| \leq M \min\big\{|u|^m,|u|^{m+1}\big\} 
\en 
with some constant $M$ independent of $t$. Now write 
$$ 
\hat g_z(t) \equiv \int_{-\infty}^{+\infty}  e^{it x}\, g_z(x)\,dx = 
\int_{-\infty}^{+\infty}  e^{it x}\, \hat V(zx)\, h(x)\,dx = 
\int_{-\infty}^{+\infty}  \hat h(t + zu)\,dV(u). 
$$ 
The assumption $\hat V(0) = \dots = \hat V^{(m)}(0) = 0$ implies that 
$\int_{-\infty}^{+\infty} (S_m \hat h)(t,u)\,dV(u) = 0$, for all $t$.  
Therefore, using (11.9), we finally get 
\bee 
|\hat g_z(t)| 
 & = &  
\bigg|\int_{-\infty}^{+\infty} \big(\hat h(t + zu) - (S_m \hat h)(t,zu)\big)\,dV(u)\bigg| \\ 
 & \leq &  
M\int_{-\infty}^{+\infty} \min\big\{|zu|^m,|zu|^{m+1}\big\}\,d|V|(u) 
 \, = \, o(z^{m+\alpha}). 
\ene 
Note that the last relation has been already discussed in the previous  
step. 
 
Thus Proposition 11.4 is proved. 
 
\vskip5mm 
{\bf Remark 11.5.} The second part of the above proof also covers  
the limit case $\alpha = 0$ of the inequality (11.5), which may be 
written as 
\be 
\bigg|\int_{-\infty}^{+\infty}  e^{it x}\, \hat V(zx)\,h(x)\,dx\bigg|\, \leq\, 
|z|^m\, \ep(z). 
\en 
More precisely, for this assertion we only need the assumptions 
$\int_{-\infty}^{+\infty} |u|^m\,|V|(du) < +\infty$, 
$\hat V(0) = \dots = \hat V^{(m)}(0) = 0$, and 
$ 
\int_{-\infty}^{+\infty} |x|^{m+1}\,|h(x)|\,dx < +\infty. 
$ 
Proposition 11.3 is irrelevant in this case.


\vskip10mm 
\section{{\bf Fourier Transforms and Fractional Derivatives}} 
\setcounter{equation}{0} 
 
\vskip2mm 
In this section we give the proof of Proposition 11.3. By its very  
definition, 
\be 
\Gamma(\alpha)\, (D_{0+}^\alpha g)(x) =  
\frac{d}{dx} \int_0^x \frac{h(t)}{(x-t)^\alpha}\, \hat V(t)\,dt, 
\en 
provided that the derivative exists (where $0 < \alpha < 1$). 
 
Given $m \geq 0$ integer, introduce the function of the real variable 
$$ 
\eta_m(t) = e^{it} - \sum_{k=0}^m \frac{(it)^k}{k!}. 
$$ 
The assumptions on $V$ imply that the first $m$ moments of the measure  
$V$ are vanishing, i.e., $\int_{-\infty}^{\infty} u^k\,dV(u) = 0$,  
for $k = 0,\dots,m$. Hence,  
$$ 
\hat V(t) = \int_{-\infty}^{\infty} \eta_m(tu)\,dV(u). 
$$ 
Respectively, changing the variable in (12.1) and applying Fubini's theorem,  
one may write 
\bee 
\Gamma(\alpha)\, (D_{0+}^\alpha g)(x)& = & 
\frac{d}{dx} \int_0^x \frac{h(t)}{(x-t)^\alpha}\,  
\bigg[\int_{-\infty}^{+\infty} \eta_m(tu)\,dV(u)\bigg]\,dt \\ 
& = & 
\frac{d}{dx} \int_{-\infty}^{+\infty} \bigg[\int_0^x h(x-t)\,  
\eta_m\big((x-t)u\big)\,\frac{dt}{t^\alpha}\bigg] \,dV(u). 
\ene 
 
We intend to move the differentiation inside the outer integral. To justify  
this step, consider the derivative with respect to the inner integral, 
$$ 
I(x,u) = \frac{d}{dx} \int_0^x  h(x-t)\, \eta_m\big((x-t)u\big)\,\frac{dt}{t^\alpha}. 
$$

\vskip5mm 
{\bf Lemma 12.1.} {\it Let $h(x)$ denote a continuously differentiable  
function on the real line, such that\, $|h(x)| + |h'(x)| = O(x^{-(2+m+\alpha)})$, 
as $|x| \rightarrow \infty$. Then, for all $u \in \R$ and $x > 0$, 
\bee 
|I(x,u)| & \leq & C (1+x)^{-(1+\alpha)}\, \min\{|u|, \, |u|^\alpha\} \hskip12mm  (m=0), \\ 
|I(x,u)| & \leq & C (1+x)^{-(1+\alpha)}\, \min\{|u|^m,|u|^{m+1}\} \ \ \ \ (m \geq 1) 
\ene 
with some constant $C$, depending on $h$ and $\alpha$, only. 
} 
 
\vskip5mm 
{\bf Proof.} Put $\xi_u(t) = h(t) \eta_m(ut)$, so that 
$$ 
I(x,u) = \frac{d}{dx} \int_0^x  \xi_u(x-t)\,\frac{dt}{t^\alpha}. 
$$ 
In this case we may interchange differentiation and integration. 
Thus, using $\frac{d}{dx}\, \xi_u(x-t) = - \frac{d}{dt}\, \xi_u(x-t)$  
together with $\eta_m(0)=0$, we can write 
$$ 
I(x,u) = - \int_0^x  \xi_u'(x-t)\,\frac{dt}{t^\alpha}. 
$$ 
 
Assume that $x \geq 1$ and split the integration into the two regions  
such that 
$$ 
I(x,u) = I_0(x,u) + I_1(x,u) = - \int_0^1  \xi_u'(x-t)\,\frac{dt}{t^\alpha} 
- \int_1^x  \xi_u'(x-t)\,\frac{dt}{t^\alpha}. 
$$ 
 
{\bf The integral} $I_1$. 
 
Integrating by parts, we have 
\be 
I_1(x,u) = \alpha \int_1^x  \xi_u(x-t)\,\frac{dt}{t^{1+\alpha}} - \xi_u(x-1). 
\en 
 
To analize this integral, we use the elementary bound 
\be 
|\eta_m(t)| \leq 4 \min\big\{|t|^m,|t|^{m+1}\big\}, \qquad t \in \R. 
\en 
Indeed, from Taylor's formula it follows that $|\eta_m(t)| \leq  
\frac{|t|^{m+1}}{(m+1)!}$. This settles (12.3) in case $|t| \leq 1$. 
In the other case $|t| \geq 1$, just write 
$$ 
|\eta_m(t)| \leq 1 + \sum_{k=0}^m \frac{|t|^k}{k!} \leq |t|^m 
\bigg(1 + \sum_{k=0}^m \frac{1}{k!}\bigg) < (1+e)\,|t|^m, 
$$ 
thus proving (12.3). 
 
This bound implies that 
\be 
|\xi_u(x-t)| \leq 4\,|h(x-t)| \min\big\{|(x-t)u|^m, |(x-t)u|^{m+1}\big\}. 
\en 
 
By the assumption on $h$, we have $|h(x-1)| \leq C x^{-(2+m+\alpha)}$, so, by (12.4), 
\be 
|\xi_u(x-1)| \leq Cx^{-(1+\alpha)} \min\big\{|u|^m, |u|^{m+1}\big\} 
\en 
with some constant $C$.  
 
In the region $1 \leq t \leq x_1 = \frac{1+x}{2}$, we use the bound 
$|h(x-t)| \leq C x^{-(2+m+\alpha)}$ with a constant independent of  
$t$ and $x$. Hence, by (12.4), in this region 
$$ 
|\xi_u(x-t)| \leq C x^{-(1+\alpha)} \min\big\{|u|^m, |u|^{m+1}\big\}, 
$$ 
and 
$$ 
\int_1^{x_1} |\xi_u(x-t)|\,\frac{dt}{t^{1+\alpha}} \leq  
C x^{-(1+\alpha)} \min\big\{|u|^m, |u|^{m+1}\big\}. 
$$ 
 
In the second region $x_1 \leq t \leq x$, just use  
$\frac{1}{t^{1+\alpha}} \leq Cx^{-(1+\alpha)}$. Then, by (12.4), 
\bee 
\int_{x_1}^x |\xi_u(x-t)|\,\frac{dt}{t^{1+\alpha}} 
 & \leq &  
4C x^{-(1+\alpha)}\, |u|^m \int_{x_1}^x |h(x-t)|\,(x-t)^m\,dt \\ 
 & = & 
4C x^{-(1+\alpha)}\, |u|^m  \int_0^{(x-1)/2} |h(t)| \, t^m\, dt \, \leq \,  
C' x^{-(1+\alpha)}\, |u|^m, 
\ene 
since the last integral is uniformly bounded. Similarly, again by (12.4), 
\bee 
\int_{x_1}^x |\xi_u(x-t)|\,\frac{dt}{t^{1+\alpha}} 
 & \leq &  
4C x^{-(1+\alpha)}\, |u|^{m+1} \int_{x_1}^x |h(x-t)|\,(x-t)^{m+1}\,dt \\ 
 & \leq & 
C' x^{-(1+\alpha)}\, |u|^{m+1}. 
\ene 
Collecting the bounds for the two regions, we get 
$$ 
\int_1^x |\xi_u(x-t)|\,\frac{dt}{t^{1+\alpha}} \leq 
C x^{-(1+\alpha)} \min\big\{|u|^m, |u|^{m+1}\big\} 
$$ 
with some constant $C$. Applying it together with (12.5) in (12.2), we arrive at 
\be 
I_1(x,u) \leq  C x^{-(1+\alpha)} \min\big\{|u|^m, |u|^{m+1}\big\}. 
\en 
 
{\bf The integral} $I_0$. 
 
Now, let us turn to the integral  
$I_0(x,u) = - \int_0^1  \xi_u'(x-t)\,\frac{dt}{t^\alpha}$.  
After differentiation and using the identity $\eta_m' = i\eta_{m-1}$  
(with the convention that $\eta_{-1}(t) = e^{it}$), one may represent  
it as $I_0 = I_{0,1} + I_{0,2}$, where 
\bee 
I_{0,1}(x,u) 
 & = &  
- iu \int_0^1 \eta_{m-1}((x-t)u)\, h(x-t)\,\frac{dt}{t^\alpha}, \\ 
I_{0,2}(x,u) 
 & = & 
- \int_0^1 \eta_m((x-t)u)\, h'(x-t)\,\frac{dt}{t^\alpha}. 
\ene 
 
By (12.3), since $x-t \leq x$, 
$$ 
|I_{0,2}(x,u)| \leq 4\int_0^1 |h'(x-t)|\,\min\big\{|xu|^m, |xu|^{m+1}\big\}\, 
\frac{dt}{t^\alpha}. 
$$ 
Using the assumption $h'(x) = O\big(x^{-(2+m+\alpha)}\big)$, we get 
\be 
|I_{0,2}(x,u)| \leq C x^{-(1+\alpha)}\, \min\big\{|u|^m, |u|^{m+1}\big\} 
\en 
with some constant $C$. 
 
As for the integral $I_{0,1}$, first rewrite it as 
$I_{0,1}(x,u) = - iu \int_0^1 h(x-t)\, d \zeta_{m-1}(t)$, where 
$$ 
\zeta_{m-1}(t) = \int_0^{t}\frac{\eta_{m-1}((x-w)u)}{w^\alpha}\, dw. 
$$ 
Integrating by parts, one may represent it as 
$I_{0,1} = I_{0,1,1} + I_{0,1,2}$, where 
\bee 
I_{0,1,1}(x,u) 
 & = &  
- iu\, h(x-1)\, \zeta_{m-1}(1), \\ 
I_{0,1,2}(x,u) 
 & = & 
- iu \int_0^1 h'(x-t)\,\zeta_{m-1}(t) dt. 
\ene 
 
{\bf Claim.} {\it For all $x \geq 1$, $t \in [0,1]$, we have  
$|\zeta_{-1}(t)| \leq C \min\{1, |u|^{\alpha - 1}\}$ with some constant~$C$,  
while in case $m \geq 1$, 
$$ 
|\zeta_{m-1}(t)| \leq C x^m \min\big\{|u|^{m-1},\,|u|^m\big\}. 
$$ 
} 
 
{\bf Proof.} If $m = 0$, whenever $u \neq 0$, 
$$ 
|\zeta_{-1}(t)| = \bigg|\int_0^t e^{-iwu}\,\frac{dw}{w^\alpha}\bigg| = 
|u|^{\alpha - 1} \bigg|\int_0^{tu} e^{-iw}\,\frac{dw}{w^\alpha}\bigg| \leq 
C\,|u|^{\alpha - 1}. 
$$ 
On the other hand, $|\zeta_{-1}(t)| \leq \int_0^1 \frac{dw}{w^\alpha}$, so 
$|\zeta_{-1}(t)| \leq C \min\{1, |u|^{\alpha - 1}\}$. 
 
\vskip2mm 
If $m \geq 1$, we just appeal to the estimate (12.4), which, for all $w \in (0,1)$  
and $x \geq 1$, yields 
\bee 
|\eta_{m-1}((x-w)u)|  
 & \leq & 
4 \min\big\{(x-w)^{m-1}\,|u|^{m-1},\, (x-w)^m\,|u|^m\big\} \\ 
 & \leq & 
4 x^m \min\big\{|u|^{m-1},\,|u|^m\big\}. 
\ene 
It immediately implies the desired estimate. 
 
\vskip5mm 
{\bf Continuation of the proof of Lemma 12.1.} 
Now, by the assumption on $h$, and using the claim in case $m=0$, we obtain  
$$ 
|I_{0,1}(x,u)| \leq |I_{0,1,1}(x,u)| + |I_{0,1,2}(x,u)| \leq  
C\, x^{-(1+\alpha)} \min\big\{|u|, |u|^{\alpha}\big\}. 
$$ 
Similarly, in case $m \geq 1$, 
\bee 
|I_{0,1,1}(x,u)| + |I_{0,1,2}(x,u)|  
 & \leq & 
C |u|\, \big(|h(x-1)| + |h'(x-1)|\big)\cdot x^m \min\big\{|u|^{m-1},\,|u|^m\big\} \\ 
 & \leq & 
C\, x^{-(1+\alpha)} \min\big\{|u|^{m},\,|u|^{m+1}\big\}. 
\ene 
Thus, 
\bee 
|I_{0,1}(x,u)| & \leq & C x^{-(1+\alpha)} \min\big\{|u|, |u|^{\alpha}\big\}\hskip10mm (m=0), \\ 
|I_{0,1}(x,u)| & \leq & C x^{-(1+\alpha)} \min\{|u|^m,|u|^{m+1}\} \ \ \ \ (m \geq 1). 
\ene 
Taking into account (12.7) and (12.6), we arrive at similar bounds  
for $I_0(x,u)$ and $I(x,u)$, 
which are equivalent forms of the desired bounds in the lemma in case $x \geq 1$. 
 
Finally, let us only note that the case $0 < x < 1$ may be treated in  
a similar manner (with simpler estimates). Thus Lemma 12.1 is proved.

\vskip5mm 
{\bf Proof of Proposition 11.3.} 
Finally, we are prepared to justify the differentiation step. Define 
$$ 
\psi(x) = \int_{-\infty}^{+\infty} 
\bigg[\int_0^x  \xi_u(x-t)\,\frac{dt}{t^\alpha}\bigg]\,dV(u), 
$$ 
where, as before, $\xi_u(t) = h(t) \eta_m(tu)$. Given 
$x>0$ and $\ep_n \rightarrow 0$ (with $\ep_n \neq 0$, $x + \ep_n > 0$), write 
$$ 
\frac{\psi(x+\ep_n) - \psi(x)}{\ep_n} = 
\int_{-\infty}^{+\infty} \bigg[\frac{1}{\ep_n}\int_x^{x+\ep_n} I(y,u)\,dy\bigg]\ dV(u). 
$$ 
By Lemma 12.1 in case $m=0$, the expression in the square brackets is bounded  
in absolute value by 
$$ 
C\, \min\{|u|,|u|^\alpha\}\,  \bigg| \frac{1}{\ep_n} \int_x^{x+\ep_n}  
(1+y)^{-(1+\alpha)}\,dy\bigg| \leq C\, \min\{|u|,|u|^\alpha\}. 
$$ 
On the right-hand side the function is integrable with respect to the measure  
$|V|$, according to the moment condition on the function $V$. 
A similar conclusion holds in case $m \geq 1$ (with appropriate 
modifications in the estimate). 
Therefore, one may apply the Lebesgue dominated convergence theorem,  
which gives 
$$ 
\psi'(x) = \lim_{n \rightarrow \infty} \frac{\psi(x+\ep_n) - \psi(x)}{\ep_n} = 
\int_{-\infty}^{+\infty} I(x,u)\,dV(u). 
$$ 
 
Thus, the fractional derivative 
$$ 
(D_{0+}^\alpha g)(x) = \frac{1}{\Gamma(\alpha)}\, \psi'(x) = \frac{1}{\Gamma(\alpha)} 
\int_{-\infty}^{+\infty} I(x,u)\,dV(u) 
$$ 
is well-defined for all $x>0$. Moreover, from Lemma 12.1 we also obtain that 
$$ 
|(D_{0+}^\alpha g)(x)| \leq C\, (1+x)^{-(1+\alpha)} 
\int_{-\infty}^{+\infty} \min\{|u|^{m_1},|u|^{m_2}\}\, d|V|(u), 
$$ 
where $m_1 = 1$, $m_2 = \alpha$ in case $m = 0$, and $m_1 = m$, $m_2 = m+1$ in case 
$m \geq 1$.  
 
This proves the first assertion and the inequality (11.3) in Proposition 11.3. 
For the second assertion, apply Proposition 11.1 and the formula (11.2)  
for the fractional integration by parts with the functions  
$f(x) = e^{-(\ep - it)\, x}$ ($\ep>0$) and use the second formula in (11.1) with  
$\lambda = \ep - it$ for the fractional derivatives of $f$. Then we obtain  
$$ 
\int_0^{+\infty} e^{-(\ep - it)\, x}\, (D_{0+}^\alpha g)(x)\,dx =  
(\ep - it)^{\alpha} \int_0^{+\infty} e^{-(\ep - it)\, x}\, g(x)\,dx. 
$$ 
Letting $\ep \rightarrow 0$ and using the integrability of both $g$ and  
$D_{0+}^\alpha g$ (due to (11.3)), we arrive in the limit at the required  
equality (11.4). Thus Proposition 11.3 is proved.


\vskip10mm 
\section{{\bf Binomial Decomposition of Convolutions}} 
\setcounter{equation}{0} 
 
\vskip2mm 
We shall now the probability densities $\widetilde \rho_n$ in Theorem 1.2.  
The next procedure is known; a related approach has been used, e.g.,  
in [S-M], [I-L] to study the central limit theorem 
with respect to the total variation distance. 
 
Let $0 < c < 1$ be a prescribed number, $m = [s]$, and $n \geq m+2$. 
 
Without loss of generality, one may assume that $n_0 = 1$, that is,  
$S_1 = X_1$ has a density, say, $\rho$ which may or may not be bounded. 
For definiteness, assume it is (essentially) unbounded, so that the integral 
$$ 
b = \int_{\{\rho(x) > M\}} \rho(x)\,dx 
$$ 
is positive for all $M>0$. We choose $M$ to be sufficiently large to satisfy,  
e.g., $0 < b < \frac{c}{2}$ which implies 
$2 n^{m+1}\,b^{n-m-1} < c^n$, for all $n \geq n_1$ large enough. 
 
Consider the decomposition 
$$ 
\rho(x) = a p(x) + b q(x), 
$$ 
where $a = 1-b$, and $p(x)$, $q(x)$ are the normalized restrictions of $\rho$  
to the sets $\{\rho(x) \leq M\}$ and $\{\rho(x) > M\}$, respectively. Hence, 
for the convolutions we have a binomial decomposition 
$$ 
\rho^{*n} = \sum_{k=0}^n C_n^k\, a^k b^{n-k}\, p^{*k} * q^{*(n-k)}. 
$$ 
Then split the above sum into the two parts to get  
$\rho^{*n}(x) = p_n(x) + q_n(x)$, where 
$$ 
p_n = \sum_{k = m+2}^n C_n^k\, a^k b^{n-k}\, p^{*k} * q^{*(n-k)}, \qquad 
q_n = \sum_{k=0}^{m+1} C_n^k\, a^k b^{n-k}\, p^{*k} * q^{*(n-k)}. 
$$ 
Note that 
$$ 
\beta_n \equiv \int_{-\infty}^{+\infty} q_n(x)\,dx =  
\sum_{k=0}^{m+1} C_n^k\, a^k b^{n-k} 
\leq n^{m+1}\,b^{n-m-1} < \frac{c^n}{2} \qquad (n \geq n_1). 
$$ 
Finally define  
$$ 
\widetilde \rho_n(x) = \frac{\sqrt{n}}{1 - \beta_n}\, p_n\big(x\sqrt{n}\big), \qquad 
\widetilde v_n(t) = \int_{-\infty}^{+\infty} e^{itx} \widetilde \rho_n(x)\,dx. 
$$ 
 
Let us recall that $\rho_n(x) = \sqrt{n}\, \rho^{*n}\big(x\sqrt{n}\big)$  
has the characteristic function $v_n(t) = v(\frac{t}{\sqrt{n}})^n$,  
where $v$ is the characteristic function of $X_1$.  
By construction, the densities $\widetilde \rho_n$ are bounded and provide 
a strong approximation for $\rho_n$. Namely, we immediately obtain: 
 
\vskip5mm 
{\bf Lemma 13.1.} {\it For all $n \geq n_1$, 
$\int_{-\infty}^{+\infty} |\widetilde \rho_n(x) - \rho_n(x)|\,dx < c^n$. 
In particular, for all $t \in \R$, 
$$ 
|\widetilde v_n(t) - v_n(t)| < c^n. 
$$ 
A similar inequality also holds for the first $m$ derivatives of  
$\widetilde v_n$ and $v_n$ with $n$ large enough. 
} 
 
\vskip5mm 
The last assertion of the lemma needs a more detailed explanation 
which we postpone to the end of the section. 
 
We will also need some integrability properties for $\widetilde v_n$ 
and their first $m$ derivatives that are due to the boundedness of the probability 
density $p(x)$ and the finiteness of the $m$-th absolute moment of $X_1$. 
 
\vskip5mm 
{\bf Lemma 13.2.} {\it Let $\E\, |X_1|^m < +\infty$, $m \geq 2$. 
There exist positive constants $A$ and $\sigma$, depending on $X_1$ and $m$, 
such that, for all $0 \leq T \leq \sqrt{n}$, 
\be 
\int_{\{|t| \geq T\}} |\widetilde v_n(t)|\,dt < A\, e^{-\sigma^2 T^2}. 
\en 
A similar bound is also true for the first $m$ derivatives of  
$\widetilde v_n$ $($with arbitrary $n \geq m+2)$. 
}

\vskip5mm 
{\bf Proof.} Let $\hat p$, $\hat q$ denote the Fourier transforms of  
$p$ and $q$, respectively. Then 
\be 
\widetilde v_n(t) =  
\frac{1}{1 - \beta_n}\sum_{k = m+2}^n C_n^k\, a^k b^{n-k}\,  
\hat p\big(\frac{t}{\sqrt{n}}\big)^k\,  
\hat q\big(\frac{t}{\sqrt{n}}\big)^{n-k}. 
\en 
 
 
By the Riemann-Lebesgue theorem,  
\be 
\sup_{|t| \geq 1} |\hat p(t)| < \gamma, \quad   
\sup_{|t| \geq 1} |\hat q(t)| < \gamma \qquad  (0 \leq \gamma < 1). 
\en 
Hence, from (13.2), for all $|t| \geq \sqrt{n}$, 
\be 
|\widetilde v_n(t)| < \frac{\gamma^{n-2}}{1 - \beta_n}\, 
\big|\hat p\big(\frac{t}{\sqrt{n}}\big)\big|^2 
\sum_{k = m+2}^n C_n^k\, a^k b^{n-k} \leq \gamma^{n-2}\,  
\big|\hat p\big(\frac{t}{\sqrt{n}}\big)\big|^2. 
\en 
In addition, by the Plancherel theorem, and using $p(x) \leq M/a$, we have 
\be 
\int_{-\infty}^{+\infty}  
\big|\hat p\big(\frac{t}{\sqrt{n}}\big)\big|^2\,dt = 2\pi \sqrt{n} \int_{-\infty}^{+\infty} p(x)^2\,dt \leq \frac{2\pi M}{a}\,\sqrt{n}. 
\en 
Therefore, integrating the inequality (13.4), we get 
\be 
\int_{\{|t| \geq \sqrt{n}\,\}} |\widetilde v_n(t)|\,dt < 
\frac{2\pi M}{a}\,\gamma^{n-2} \sqrt{n}. 
\en 
 
On the other hand, since both $p$ and $q$ represent the densities of 
probability distributions with finite second moments,  
their characteristic functions near zero satisfy 
\be 
|\hat p(t)| \leq e^{-\sigma^2 t^2}, \quad |\hat q(t)| \leq e^{-\sigma^2 t^2}  
\qquad (|t| \leq 1) 
\en 
with some constant $\sigma>0$. Hence, for $|t| \leq \sqrt{n}$, (13.2) gives  
the estimate $|\widetilde v_n(t)| \leq e^{-\sigma^2 t^2}$ and 
$$ 
\int_{\{T \leq |t| \leq \sqrt{n}\,\}} |\widetilde v_n(t)|\,dt \, \leq \,  
\int_{T \leq |t| \leq \sqrt{n}}\, e^{-\sigma^2 t^2}\,dt \, 
< \, \frac{1}{\sigma}\,e^{-\sigma^2 T^2}. 
$$ 
Together with (13.6) the latter yields 
$$ 
\int_{\{|t| \geq T\,\}} |\widetilde v_n(t)|\,dt < 
\frac{1}{\sigma}\,e^{-\sigma^2 T^2} + \frac{2\pi M}{a}\,\gamma^{n-2} \sqrt{n}. 
$$ 
Finally, since for the values $0 \leq T \leq \sqrt{n}$ one always has  
$\gamma^{n-2} \sqrt{n} \leq A_1\, e^{-\sigma_1^2 T^2}$ with some constants  
$A_1$ and $\sigma_1>0$, the desired bound (13.1) easily follows. 
 
\vskip2mm 
As for the derivatives, a bound of this type can be proved by similar arguments,  
so let us restrict ourselves to the basic case of the $m$-th derivative  
(needed for the proof of Theorem~1.2). 
 
The condition $\E |X_1|^m = \int |x|^m\,\rho(x)\,dx < +\infty$ implies 
a similar property for the densities $p(x)$ and $q(x)$. Hence,  
$\hat p(t)$ and $\hat q(t)$ have continuous derivatives up to order $m$,  
bounded in absolute value by some common constant. 
 
In view of (13.2), $\widetilde v_n^{(m)}(t)$ represents a linear  
combination of the terms 
\be 
\frac{d^m}{dt^m}\, 
\bigg[\hat p\big(\frac{t}{\sqrt{n}}\big)^k\, \hat q\big(\frac{t}{\sqrt{n}}\big)^{n-k}\bigg] 
= \sum_{r=0}^m C_m^r \, 
\frac{d^r}{dt^r}\, 
\bigg[\hat p\big(\frac{t}{\sqrt{n}}\big)^k\bigg]\, 
 \frac{d^{m-r}}{dt^{m-r}}\,\bigg[\hat q\big(\frac{t}{\sqrt{n}}\big)^{n-k}\bigg] 
\en 
with integers $m+2 \leq k \leq n$. Given $0 \leq r \leq m$, by the chain  
rule (2.3), the $r$-th derivative of the function  
$\hat p(\frac{t}{\sqrt{n}})^k$ 
represents a linear combination of the terms 
\be 
n^{-r/2}\,\hat p\big(\frac{t}{\sqrt{n}}\big)^{k - (k_1 + \dots + k_r)}\,  
\hat p'\big(\frac{t}{\sqrt{n}}\big)^{k_1} \dots  
\hat p^{(r)}\big(\frac{t}{\sqrt{n}}\big)^{k_r} 
\en 
over all integer tuples $(k_1,\dots,k_r)$, such that  
$k_1 + 2k_2 + \dots + r k_r = r$ and $k_1 + \dots + k_r \leq k$  
($k_j \geq 0$). Moreover, the coefficients in that linear combination  
do not depend on $n$, and the total number of such terms 
is bounded by a quantity that depends on $m$, only. 
 
Using $k_1 + \dots + k_r \leq r \leq m$ and the boundedness of the derivatives, 
the absolute value of the expression (13.9), as well as the sum of all  
such terms is bounded by $|\hat p(\frac{t}{\sqrt{n}})|^{k - m}$  
up to a constant factor. 
 
Similarly, the $(m-r)$-th derivative of the function $\hat q(\frac{t}{\sqrt{n}})^{n-k}$ 
represents a linear combination of the terms 
\be 
n^{-(m-r)/2}\, 
\hat q\big(\frac{t}{\sqrt{n}}\big)^{(n-k) - (k_1 + \dots + k_{m-r})}\,  
\hat q'\big(\frac{t}{\sqrt{n}}\big)^{k_1} \dots  
\hat q^{(m-r)}\big(\frac{t}{\sqrt{n}}\big)^{k_{m-r}} 
\en 
over all integer tuples $(k_1,\dots,k_{m-r})$, such that  
$k_1 + 2k_2 + \dots + (m-r) k_{m-r} = m-r$ 
and $k_1 + \dots + k_{m-r} \leq n-k$ ($k_j \geq 0$). 
Again, the coefficients in the linear combination do not depend on $n$, and  
the total number of such terms is bounded by a quantity depending on $m$, only. 
Since $k_1 + \dots + k_{m-r} \leq \min(n-k,m-r) \leq \min(n-k,m)$, the absolute  
value of the expression (13.10), and the sum of all such terms is bounded by  
$|\hat q(\frac{t}{\sqrt{n}})|^{(n-k) - \min(n-k,m)}$ up to a constant factor. 
 
Thus, (13.8) is bounded in absolute value by 
\be 
C\, \big|\hat p(\frac{t}{\sqrt{n}})\big|^{k - m}\ 
\big|\hat q(\frac{t}{\sqrt{n}})\big|^{(n-k) - \min(n-k,m)} 
\en 
with some constant $C$, depending on $X_1$ and $m$, only. It then follows from  
(13.2) that  
\be 
\big|\widetilde v_n^{(m)}(t)\big| \, \leq \, \frac{C}{1-\beta_n}  
\sum_{k = m+2}^n C_n^k\, a^k b^{n-k}\,  
\big|\hat p(\frac{t}{\sqrt{n}})\big|^{k - m}\ 
\big|\hat q(\frac{t}{\sqrt{n}})\big|^{(n-k) - \min(n-k,m)}. 
\en 
 
Now, like in the previous step, using (13.3), for all  
$|t| \geq \sqrt{n}$, we get 
\bee 
\big|\widetilde v_n^{(m)}(t)\big| 
 & \leq &  
\frac{C}{1-\beta_n}\ \big|\hat p\big(\frac{t}{\sqrt{n}}\big)\big|^2 
\sum_{k = m+2}^n C_n^k\, a^k b^{n-k} \, \gamma^{(n-m-2) - \min(n-k,m)} \\ 
 & \leq &  
C\, \gamma^{n-2m-2}\, \big|\hat p\big(\frac{t}{\sqrt{n}}\big)\big|^2. 
\ene 
Integrating this inequality with the help of (13.5), we obtain that 
\be 
\int_{\{|t| \geq \sqrt{n}\,\}} |\widetilde v_n^{(m)}(t)|\,dt < 
\frac{2\pi MC}{a}\,\gamma^{n-2m-2} \sqrt{n}. 
\en 
In addition, using (13.7), for $|t| \leq \sqrt{n}$, the product (13.11) is  
bounded by $C\,e^{-d\sigma^2 t^2}$, where 
$$ 
d = \frac{1}{n}\,\big((n-m) - \min(n-k,m)\big) \geq  
d' = \frac{1}{n}\,\big((n-m) - \min(n-m-2,m)\big). 
$$ 
If $n \geq 2m+2$, then $d' = \frac{n-2m}{n} \geq \frac{1}{m+1}$.  
In the other case $m+2 \leq n < 2m+2$, we also have  
$d' = \frac{2}{n} \geq \frac{1}{m+1}$. Hence, (13.11) is bounded by 
$C\,e^{-\sigma^2 t^2/(m+1)}$, and we derive from (13.12) 
$$ 
\big|\widetilde v_n^{(m)}(t)\big| \, \leq \, \frac{C}{1-\beta_n}\,  
\sum_{k = m+2}^n C_n^k\, a^k b^{n-k}\,  
e^{-\sigma^2 t^2/(m+1)} \, \leq \, C\, e^{-\sigma^2 t^2/(m+1)}. 
$$ 
It remains to integrate this inequality to get 
$$ 
\int_{\{T \leq |t| \leq \sqrt{n}\,\}} |\widetilde v_n^{(m)}(t)|\,dt < 
C \int_{T \leq |t| \leq  
\sqrt{n}} e^{-\sigma^2 t^2/(m+1)}\,dt < 
\frac{C\sqrt{m+1}}{\sigma}\,e^{-\sigma^2 T^2/(m+1)}. 
$$ 
Together with (13.13), it yields the desired estimate (13.1).  
Thus Lemma 13.2 is proved.

\vskip5mm 
{\bf Remark 13.3.} If the density $\rho$ is bounded, the decomposition  
procedure is not needed, and then Lemma 13.2 should read as follows.  
Let $\E\, |X|^m < +\infty$, $m \geq 2$, for a random variable having  
a bounded density. There exist constants $A$ and $\sigma>0$, such that  
for all $n \geq 2$, 
\be 
\int_{\{|t| \geq T\}} \big|v\big(\frac{t}{\sqrt{n}}\big)\big|^n\,dt <  
A\, e^{-\sigma^2 T^2}, \qquad 0 \leq T \leq \sqrt{n}, 
\en 
where $v$ is the characteristic function of $X$. 
A similar bound holds as well for the first $m$ derivatives of  
$v(\frac{t}{\sqrt{n}})^n$ with arbitrary $n \geq m+2$.

\vskip5mm 
{\bf Proof of Lemma 13.1.} By the construction, for all $n \geq n_1$, 
$$ 
\int_{-\infty}^{+\infty} |\widetilde \rho_n(x) - \rho_n(x)|\,dx \leq  
2\beta_n < c^n, 
$$ 
so $|\widetilde v_n(t) - v_n(t)| < c^n$, as well. In order to extend  
this inequality to the derivatives, recall the representation (13.2)  
to write 
\be 
\widetilde v_n(t) - v_n(t) =  
\frac{\beta_n}{1 - \beta_n}\, \Sigma_1 - \Sigma_2, 
\en 
where 
\bee 
\Sigma_1 & = &  
\sum_{k = m+2}^n C_n^k\, a^k b^{n-k}\,  
\hat p\big(\frac{t}{\sqrt{n}}\big)^k\, \hat q\big(\frac{t}{\sqrt{n}}\big)^{n-k}, \\ 
\Sigma_2 & = &  
\sum_{k = 0}^{m+1}\, C_n^k\, a^k b^{n-k}\,  
\hat p\big(\frac{t}{\sqrt{n}}\big)^k\, \hat q\big(\frac{t}{\sqrt{n}}\big)^{n-k}. 
\ene 
As before, we will only consider the case of the $m$-th derivative. 
 
It has been shown in the proof of Lemma 13.2 that, given  
$m+2 \leq k \leq m$, the function 
$\hat p(\frac{t}{\sqrt{n}})^k\, \hat q\big(\frac{t}{\sqrt{n}})^{n-k}$  
has the $m$-th derivative, bounded in absolute value by the expression  
(13.11). So, it is bounded by a constant $C$, depending on $X_1$ and $m$,  
only. In the general case, including the values $0 \leq k \leq m+1$,  
(13.11) should be replaced with 
$$ 
C\, \big|\hat p(\frac{t}{\sqrt{n}})\big|^{\max(k - m,0)}\ 
\big|\hat q(\frac{t}{\sqrt{n}})\big|^{(n-k) - \min(n-k,m)}, 
$$ 
which is also bounded by $C$. Therefore, from (13.15), 
\bee 
\big|\widetilde v_n^{(m)}(t) - v_n^{(m)}(t)\big| 
 & \leq &  
\frac{C\beta_n}{1 - \beta_n}\, \sum_{k = m+2}^n C_n^k\, a^k b^{n-k} + 
C \sum_{k = 0}^{m+1}\, C_n^k\, a^k b^{n-k} \\ 
 & = & 2C\beta_n \, < \, c^n, 
\ene 
where the last inequality holds true for all $n$ starting with a certain  
$n_1$.  
 
Thus, Lemma 13.1 is proved.


\vskip10mm 
\section{{\bf Proof of Theorems 1.1 and 1.2}} 
\setcounter{equation}{0} 
 
\vskip2mm 
We are prepared to make the last step in the proof of Theorems 1.1 and 1.2. 
Recall that $s \geq 2$, $m = [s]$, and put $\alpha= s-m$. 
 
Let $v(t)$ be the characteristic function of $X_1$ and 
$v_n(t) =  v(\frac{t}{\sqrt{n}})^n$ be the characteristic function of $S_n$. 
We will assume that all $S_n$ have densities $\rho_n$ 
(since only minor modifications have to be done in the more general case, where 
$S_n$ have densities for all $n$ large enough). 
 
If $\rho_{n_0}$ and therefore all $\rho_n$ with $n \geq n_0$ are (essentially)  
bounded for some $n_0$, then there is no need to use the binomial decomposition  
of the previous section, and we put $\widetilde \rho_n = \rho_n$. 
This case corresponds to Theorem 1.1. 
Otherwise, if $\rho_n$ are unbounded for all $n \geq 1$, 
then the binomial decomposition is applied to $\rho = \rho_1$, and 
we obtain the modified densities $\widetilde \rho_n$ together 
with the associated characteristic functions $\widetilde v_n$, which we considered  
in the previous section. Thus, the requirement $c)$ in Theorem 1.2 is met.

\vskip2mm 
{\bf The inversion formula.} 
 
The characteristic functions $\widetilde v_n$ have continuous, bounded derivatives  
up to order $m$, that are integrable according to the inequality (13.1) of Lemma 13.2,  
or (13.14) of Remark 13.3. Hence, by the inversion formula, 
$$ 
(ix)^p\ \widetilde\rho_n(x) = \frac{1}{2\pi} \int_{-\infty}^{+\infty}  
e^{-itx}\, \widetilde v_n^{(p)}(t)\,dt, \qquad p = 0,1,\dots,m. 
$$ 
 
By the construction, the approximating functions 
$\varphi_m(x) = \varphi(x) + \sum_{k=1}^{m-2} q_k(x)\,n^{-k/2}$, which appear in 
the relation (1.3), have the integrable Fourier transform 
$$ 
u_m(t) = e^{-t^2/2}\, \bigg(1 + \sum_{k=1}^{m-2} P_k(it)\, n^{-k/2}\bigg). 
$$ 
Consequently, for all $p = 0,1,\dots,m$, 
\be 
(ix)^p\, \big(\widetilde \rho_n(x) - \varphi_m(x)\big) =  
\frac{1}{2\pi} \int_{-\infty}^{+\infty}  
e^{-itx} \big(\widetilde v_n^{(p)}(t) - u_m^{(p)}(t)\big)\,dt. 
\en 
Our task is thus to give proper upper bounds on the absolute value  
of these integrals in the particular cases $p = 0$ and $p = m$.  
 
What is rather standard, one should split the integration into two regions. 
Given $T_n \rightarrow +\infty$, $0 \leq T_n \leq \sqrt{n}$ (to be specified 
later on), let 
$$ 
I_{n,p} = \int_{|t| \leq T_n} 
e^{-itx} \big(\widetilde v_n^{(p)}(t) - u_m^{(p)}(t)\big)\,dt, \qquad  
J_{n,p} = \int_{|t| \geq T_n}  
e^{-itx} \big(\widetilde v_n^{(p)}(t) - u_m^{(p)}(t)\big)\,dt. 
$$ 
It should be clear that  
$$ 
\int_{|t| \geq T_n} \big|u_m^{(p)}(t)\big|\,dt \leq A e^{-\sigma^2 T_n^2} 
$$ 
with some positive constants $A$ and $\sigma$, depending on $m$. 
By Lemma 13.2 and Remark 13.3, we have a similar bound for $\widetilde v_n^{(p)}(t)$,  
whenever $n \geq m+2$, so 
\be 
|J_{n,p}| \leq A\,e^{-\sigma^2 T_n^2}. 
\en

{\bf The integral} $I_{n,p}$. 
 
To treat this integral, we subtract and add $v_n^{(p)}(t)$ inside the integrand and  
apply Lemma 13.1 (the second part). Then it gives (using $T_n \leq \sqrt{n}$) 
$$ 
|I_{n,p}| \leq |I_{n,p}'| + c^n \sqrt{n}, 
$$ 
where $0<c<1$ is the prescribed parameter in Theorem 1.2, and 
$$ 
I_{n,p}' = \int_{|t| \leq T_n} 
e^{-itx} \big(v_n^{(p)}(t) - u_m^{(p)}(t)\big)\,dt. 
$$ 
Using (14.1)-(14.2), we obtain that, for all $x$, 
\be 
|x|^p\, \big|\widetilde \rho_n(x) - \varphi_m(x)\big| \leq  
\frac{1}{2\pi}\, |I_{n,p}'| + A\,e^{-\sigma^2 T_n^2} + c^n \sqrt{n}, 
\en 
up some positive constants $A$ and $\sigma$. 
 
\vskip5mm 
{\bf Proof of (1.4) in case $|x| \geq 1$}. 
 
Note that for $|x| \leq 1$, the relation (1.4) follows from (1.3). As for the  
values $|x| \geq 1$, only the value $p=m$ is of interest in (14.3). 
 
The integral $I_{n,p}'$ can be treated with the help of Theorem 1.3.  
It gives that in the interval $|t| \leq c_1 n^{1/6}$ with some constant  
$0 < c_1 \leq 1$ we have 
$$ 
v_n^{(p)}(t) - u_m^{(p)}(t) = n\, \frac{d^p}{dt^p} 
\left[\bigg(v\big(\frac{t}{\sqrt{n}}\big) - e_m\big(\frac{t}{\sqrt{n}}\big)\bigg)\,  
e^{-t^2/2}\right] + r_n, 
$$ 
where the remainder satisfies 
$$ 
|r_n| \leq e^{-t^2/4}\, \bigg(\frac{C}{n^{(m-1)/2}} + \frac{\ep_n}{n^{s-2}}\bigg). 
$$ 
Here $C$, $c$ and $\ep_n$ are some positive constants, such that $\ep_n \rightarrow 0$, 
as $n \rightarrow \infty$. Hence, assuming that $T_n \leq c_1 n^{1/6}$ and noting 
that $c^n \sqrt{n}$ will be absorbed by other remainder terms, we get that 
\be 
|x|^p\, \big|\widetilde \rho_n(x) - \varphi_m(x)\big| \leq  
\frac{n}{2\pi}\, |I_{n,p}''| + A\,e^{-\sigma^2 T_n^2} + 
\frac{C}{n^{(m-1)/2}} + \frac{\ep_n}{n^{s-2}}, 
\en 
where 
\be 
I_{n,p}'' = \int_{|t| \leq T_n} e^{-itx} \frac{d^p}{dt^p} 
\left[\bigg(v\big(\frac{t}{\sqrt{n}}\big) - e_m\big(\frac{t}{\sqrt{n}}\big)\bigg)\,  
e^{-t^2/2}\right] dt. 
\en 
 
Now, one can differentiate inside the last integral, which will lead to the terms, 
containing $e^{-t^2/2}$ up to polynomial factors (due to the property that $v$ 
has $m$ bounded derivatives). Hence, integration in (14.5) may be extended to 
the whole real line at the expense of an error not exceeding $C e^{-T_n^2/4}$. 
Hence, (14.4) may be replaced with 
\be 
|x|^p\, \big|\widetilde \rho_n(x) - \varphi_m(x)\big| \leq  
\frac{n}{2\pi}\, |I_{n,p}'''| + A\,e^{-\sigma^2 T_n^2} + 
\frac{C}{n^{(m-1)/2}} + \frac{\ep_n}{n^{s-2}}, 
\en 
where 
$$ 
I_{n,p}''' = \int_{-\infty}^{+\infty} e^{-itx} \frac{d^p}{dt^p} 
\left[\bigg(v\big(\frac{t}{\sqrt{n}}\big) - e_m\big(\frac{t}{\sqrt{n}}\big)\bigg)\,  
e^{-t^2/2}\right] dt. 
$$ 
Letting $p=m$ and 
$w(t) = v(t) - e_m(t)$, and performing differentiation, rewrite the above integral as 
\be 
I_{n,m}'''\, =\, \sum_{k=0}^{m}\, \frac{C_m^k}{n^{k/2}}\, \int_{-\infty}^{+\infty}  
e^{-itx}\, w^{(k)}\big(\frac{t}{\sqrt{n}}\big)\, h_{m-k}(t)\, dt, 
\en 
where $h_{m-k}(t) = (e^{-t^2/2})^{(m-k)} = (-1)^{m-k} H_{m-k}(t)\, e^{-t^2/2}$. 
 
Recall that $w(t) = \hat V(t)$ represents the Fourier transform of a finite  
signed measure, $V$, such that  
$\int_{-\infty}^{+\infty} |u|^{m+\alpha}\, |V|(du) < +\infty$ 
(where $|V|$ denotes the variation of $V$, treated as a positive finite measure). 
In addition, the first $m$ derivatives of $w$ are vanishing 
(cf. Section 7 and Proposition 7.4). Hence, 
$w^{(k)}(t) = \hat V_k(t)$ represents the Fourier transform of a finite  
signed measure $V_k$, such that  
$\int_{-\infty}^{+\infty} |u|^{(m-k) + \alpha}\, |V_k|(du) < +\infty$, 
and also the first $m-k$ derivatives of $w_k$ are vanishing. 
Therefore, we are in position to apply Proposition 11.4 to the functions 
$\hat V_k(t)$ in place of $V$, $h_{m-k}$ in place of $h$, and 
with $m-k$ in place of the parameter $m$. Choosing 
$z = 1/\sqrt{n}$, the inequalities (11.5) and (11.10) (cf. Remark 11.5) give 
$$ 
\bigg|\int_{-\infty}^{+\infty}  e^{-it x}\,  
w^{(k)}\big(\frac{t}{\sqrt{n}}\big)\, h_{m-k}(t)\,dt\bigg|\, \leq\, 
\frac{\ep_n}{n^{(m-k+\alpha)/2}}\, (1 + |x|)^{-\alpha} 
$$ 
with some sequence $\ep_n \rightarrow 0$, as $n \rightarrow \infty$. 
Applying this bound in (14.7), we obtain 
$$ 
|I_{n,m}'''| \leq \frac{\ep_n}{n^{s/2}}\, (1 + |x|)^{-\alpha}, 
$$ 
and then (14.6) with $p=m$ yields 
\be 
|x|^s\, \big|\widetilde \rho_n(x) - \varphi_n(x)\big| \leq  
\frac{\ep_n}{n^{(s-2)/2}} + |x|^\alpha \bigg(A\,e^{-\sigma^2 T_n^2} + 
\frac{C}{n^{(m-1)/2}} + \frac{\ep_n}{n^{s-2}}\bigg). 
\en 
 
It remains to involve information about the possible growth of $T_n$. But, 
as we have seen, one could choose $T_n$ to be of order $n^{-1/6}$ regardless of $s$. 
With this choice (14.8) leads to the announced inequality (1.4) of Theorem 1.2.

\vskip5mm 
{\bf Proof of (1.3)}. 
 
Let us return to (14.3). The integral $I_{n,p}'$ can also be estimated by virtue  
of Proposition 9.1 and Propositions 5.1--5.2 (cf. Corollary (9.2)).  
Namely, they give that 
\be 
v_n^{(p)}(t) - u_m^{(p)}(t) = o(n^{-(s-2)/2})\, e^{-t^2/4}, \qquad p = 0,1\dots,m, 
\en 
uniformly over all $t$ in the intervals $|t| \leq T_n$, where $T_n$ are of order  
$n^{-1/6}$ in case $s \geq 3$ and of order $n^{-(s-2)/(2s)}$ in case $2 < s < 3$.  
If $s = 2$, we may only have $T_n \rightarrow +\infty$. Clearly, in all cases 
$I_{n,p}' = o(n^{-(s-2)/2})$, and (14.3) yields 
\be 
|x|^p\ |\widetilde \rho_n(x) - \varphi_m(x)| \leq  
o\big(n^{-(s-2)/2}\big) + A\,e^{-\sigma^2 T_n^2}. 
\en 
It remains to apply this inequality with $p=0$ and $p=m$. 
 
Theorems 1.1-1.2 are thus proved. 
 
\vskip5mm 
{\bf Remark.} If $\E\,|X_1|^{m+1} < +\infty$, $m \geq 2$, but $\varphi_m$  
are constructed with the help of the same cumulants $\gamma_3,\dots,\gamma_m$  
(like in the case $m \leq s < m+1$), the relation (1.3) for both Theorem 1.1 and 
Theorem 1.2 may be sharpened. Indeed, by Proposition 9.1 (second part), (14.9)  
should be replaced with a stronger relation 
$$ 
v_n^{(p)}(t) - u_m^{(p)}(t) = O(n^{-(m-1)/2})\, e^{-t^2/4}, \qquad p = 0,1\dots,m, 
$$ 
which holds uniformly in the intervals $|t| \leq c_1 n^{1/6}$. Respectively, 
it provides a stronger variant of (14.10), namely, 
$$ 
(1 + |x|^m)\ |\widetilde \rho_n(x) - \varphi_m(x)| = 
O\big(n^{-(m-1)/2}\big). 
$$

\end{document}